\documentclass[twocolumn]{autart}

\usepackage{amsmath,amssymb,mathtools}
\usepackage{graphicx,color,subfig}
\usepackage{cite}

\let\theoremstyle\relax
\usepackage{amsthm}

\usepackage{tikz,calc,bbding}
\usetikzlibrary{shapes, arrows, decorations, decorations.pathreplacing, decorations.pathmorphing, decorations.markings, fit, matrix}
\usepackage[font=small,labelfont=bf]{caption}
\usepackage{arydshln}

\graphicspath{{figures/}}

\newtheorem{theorem}{Theorem}[section]
\newtheorem{proposition}[theorem]{Proposition}
\newtheorem{lemma}[theorem]{Lemma}

\newtheorem{definition}[theorem]{Definition}
\newtheorem{example}[theorem]{Example}
\newtheorem{problem}{Problem}

\newtheorem{assumption}[theorem]{Assumption}

\newcommand{\longthmtitle}[1]{\mbox{}{\bf \textit{(#1).}}}

\newcommand{\real}{\ensuremath{\mathbb{R}}}
\newcommand{\realpos}{\ensuremath{\mathbb{R}_{>0}}}

\newcommand{\intpos}{{\mathbb{N}}}

\DeclareMathOperator*{\dlim}{dlim}

\newcommand{\Ac}{\mathcal{A}}
\newcommand{\Bc}{\mathcal{B}}

\newcommand{\Dc}{\mathcal{D}}

\newcommand{\Nc}{\mathcal{N}}

\newcommand{\Rc}{\mathcal{R}}

\newcommand\rbf{\mathbf{r}}

\newcommand\wbf{\mathbf{w}}
\newcommand\xbf{\mathbf{x}}
\newcommand\ybf{\mathbf{y}}
\newcommand\zbf{\mathbf{z}}
\newcommand\Abf{\mathbf{A}}
\newcommand\Bbf{\mathbf{B}}

\newcommand\cov{{\rm Cov}}
\newcommand\var{{\rm Var}}

\newcommand\Sigmab{\boldsymbol{\Sigma}}

\newcommand\xib{\boldsymbol{\xi}}
\newcommand\zetab{\boldsymbol{\zeta}}
\newcommand\upsilonb{\boldsymbol{\upsilon}}
\newcommand\omegab{\boldsymbol{\omega}}
\newcommand\mub{\boldsymbol{\mu}}

\newcommand\thetab{\boldsymbol{\theta}}

\newcommand\etab{{\boldsymbol{\eta}}}

\newcommand{\zeros}{\mathbf{0}}

\newcommand{\E}{\mathbb{E}}

\newcommand\li{{\rm li}}
\newcommand\Ei{{\rm Ei}}

\usepackage{enumitem}
\renewcommand{\theenumi}{(\roman{enumi})}

\renewcommand\P{\mathbb{P}}

\renewcommand{\footnoterule}{%
  \hspace{3pt} \hrule width 0.4\textwidth height 0.5pt
  \kern 2pt
}

\newcommand{\oprocendsymbol}{\hbox{$\square$}}
\newcommand{\oprocend}{\relax\ifmmode\else\unskip\hfill\fi\oprocendsymbol}
\def\eqoprocend{\tag*{$\square$}}

\makeatletter
\newcommand\inlineeqno{\stepcounter{equation}\ \@eqnnum}
\makeatother

\makeatletter
\newcommand*{\inlineequation}[2][]{%
  \begingroup
    \refstepcounter{equation}%
    \ifx\\#1\\%
    \else
      \label{#1}%
    \fi
    \relpenalty=10000 %
    \binoppenalty=10000 %
    \ensuremath{%
      #2%
    }%
    \hfill \@eqnnum
  \endgroup
}
\makeatother

\begin{document}

\begin{frontmatter}
\runtitle{Emergence of Linear Behavior via Spatial Averaging}  

\title{On the Emergence of Linear Behavior in Large-Scale Dynamical Systems via Spatial Averaging} 

\thanks[footnoteinfo]{A preliminary version of this work appeared at the 61st IEEE Conference on Decision and Control~\cite{ahmed2022linearizing}. Corresponding author Erfan Nozari.}

\author[sabbir]{Sabbir Ahmed}\ead{sahmed9@binghamton.edu},    
\author[fareed]{Hafiz Fareed Ahmed}\ead{hahme031@ucr.edu},               
\author[erfan]{Erfan Nozari}\ead{erfan.nozari@ucr.edu}  

\address[sabbir]{Department of Computer Science, State University of New York (SUNY) Binghamton, New York, NY}  
\address[fareed]{Department of Mechanical Engineering, University of California Riverside, Riverside, CA}             
\address[erfan]{Departments of Mechanical Engineering, Electrical and Computer Engineering, and Bioengineering, University of California Riverside, Riverside, CA}        


\begin{abstract}                          
Various natural and engineered systems, from urban traffic flow to the human brain, can be described by large-scale networked dynamical systems. These systems are similar in being comprised of a large number of microscopic subsystems, each with complex nonlinear dynamics and interactions, that collectively give rise to different forms of macroscopic dynamics. Despite significant research, why and how various forms of macroscopic dynamics emerge from underlying micro-dynamics remains largely unknown. In this work we focus on linearity as one of the most fundamental aspects of system dynamics. By extending the theory of mixing sequences, we show that \textit{in a broad class of autonomous nonlinear networked systems, the dynamics of the average of all subsystems' states becomes asymptotically linear as the number of subsystems grows to infinity, provided that, in addition to technical assumptions, pairwise correlations between subsystems decay to 0 as their pairwise distance grows to infinity}. We prove this result when the latter distance is between subsystems' linear indices or spatial locations, and provide extensions to linear time-invariant (LTI) limit dynamics, finite-sample analysis of rates of convergence, and networks of spatially-embedded subsystems with random locations. To our knowledge, this work is the first rigorous analysis of macroscopic linearity in large-scale heterogeneous networked dynamical systems, and provides a solid foundation for further theoretical and empirical analyses in various domains of science and engineering.
\end{abstract}

\end{frontmatter}


%
\section{Introduction}\label{sec:intro}

The spatial scale of analysis is a decisive factor in studying large-scale systems, from engineered systems %
to natural phenomena, with vast implications for modeling, system analysis, and control%
~\cite{breakspear2017dynamic,burghout2005hybrid,marchiori2015micro,lampropoulos2010methodology,jackson2017economic,banos2015coupling}. In many large-scale systems, individual microscopic subsystems exhibit complex, nonlinear dynamics that give rise, often in ways that we do not fully understand, to qualitatively distinct %
emergent dynamics at the macroscale~\cite{social_network, power_system, neural_pop, neural_pop2, neural_pop3, disease_dyn}. Therefore, in principle, macroscopic dynamics can be studied in hyper-dimensional models with microscopic resolution~\cite{markram2006blue,yao2023whole}. However, the sheer complexity of this approach often makes it infeasible,
leading many researchers to use linear models \textit{as approximations or local linearizations} %
for studying complex networks~\cite{liu2011controllability,gu2015controllability,hahn2023towards}.

In a recent empirical study~\cite{nozari2023macroscopic}, on the other hand, we found linear models to be surprisingly effective--in fact, seemingly optimal--for capturing the macroscopic dynamics of large-scale brain networks across hundreds of subjects and different data modalities.
Remarkably, this linearity appears to be in great contrast to the well-established theoretical and empirical nonlinearity of individual neuron dynamics~\cite{neuron_nonlin, neuron_nonlin1}.
Nevertheless, this observation has since been replicated by others~\cite{sani2024dissociative} as well as various simulations we have carried~\cite{ahmed2022linearizing}, and has remained begging for theoretical understanding and a rigorous mathematical framework that explains precisely why and when it holds. The goal of this work is to provide such a theoretical foundation for general stochastic nonlinear systems. %

 \subsubsection*{Literature Review}

Our initial study~\cite{nozari2023macroscopic} %
further aligns with other recent observations in large-scale biological~\cite{sani2024dissociative,yang2021modelling,sani2021modeling,schulz2019deep,liu2010linear} and artificial~\cite{luo2021phase,lee2019wide,liu2020linearity} neural networks.
In the data-driven modeling of large-scale brain dynamics linked to electrical stimulation~\cite{yang2021modelling} or motor behavior~\cite{sani2021modeling}, nonlinear models have repeatedly failed to outperform linear models, challenging the role and extent of nonlinearity in large-scale brain dynamics. More recently, it was found in~\cite{sani2024dissociative} that nonlinear readouts can indeed improve the accuracy of predicting behavioral outcomes from neural recordings, but linear models are as accurate as nonlinear ones in predictive modeling of neural data itself. Similarly, in structural and functional brain scans, nonlinear models have performed on par with linear ones in tasks such as age or sex prediction~\cite{schulz2019deep}, and correlations between visual stimuli and macroscopic brain recordings were observed to be linear in~\cite{liu2010linear}. %
Recent studies also support linear models in resting-state fMRI analysis, showing strong performance in neuropsychiatric case-control comparisons~\cite{bryant2024extracting}, and that resting-state dynamics are explained well by stationary, linear properties of the data~\cite{hosaka2024linear}.
In the realm of artificial neural networks (ANNs) with nonlinear activation functions, the linearity of gradient descent dynamics concerning network parameters has been demonstrated in the limit of infinite network width~\cite{lee2019wide,luo2021phase}. Earlier studies~\cite{jacot2018neural, liu2020linearity} also established asymptotic linearity of ANNs with respect to their parameters. Collectively, this growing body of evidence strongly suggests that, counter-intuitively, the collective dynamics of millions to billions of highly-nonlinear microscopic elements can in fact be \textit{less nonlinear} than each individual subsystem. %

In~\cite{nozari2023macroscopic}, we demonstrated through simulations that spatial averaging can be a key factor in explaining emergent macroscopic linearity. However, this was only supported by simulations of two neuron models. %
In the preliminary version of this work~\cite{ahmed2022linearizing}, we extended our simulations to more general and complex forms of microscopic nonlinear dynamics, and presented the first mathematical proof supporting the conjecture that spatial averaging induces macroscopic linearity in nonlinear systems. However, this result was shown under the strong assumption of independent and identically-distributed (i.i.d.) subsystems, in particular ruling out all forms of network interactions. %
The present work extends the methodology of~\cite{ahmed2022linearizing} to systems with heterogeneous subsystems and different forms of spatially correlated dynamics. %
Our approach leverages tools from probability theory~\cite{durrett2019probability,gut2005probability,functional_clt,jensen2000}, particularly the central limit theorem, and incorporates concepts from the theory of strong mixing sequences~\cite{bradley1986basic}, which are critical for addressing the complexities introduced by spatial correlations. %

\subsubsection*{Statement of Contributions}
The main contribution of this work is the introduction of a general theoretical framework that explains the linearizing effects of spatial averaging in populations of nonlinear dynamical systems. This can be broken into four sets of technical contributions, as follows.
\textit{First}, we extend the well-established theory of mixing stochastic processes in a number of ways, including (i) generalizing the central limit theorem (CLT) for $\rho$-mixing sequences of random variables to multivariate processes, (ii) defining and theoretically characterizing the notion of a residual factor for $\rho$-mixing sequences, and (iii) extending the definition and properties of $\rho$-mixing sequences, including the CLT, to $\rho^*$-mixing for spatially-embedded sequences with fixed as well as random Euclidean locations.
Building on these, \textit{our second and main contribution} consists of the mathematical formulation and proof of the asymptotic linearity of averaged state dynamics under spatial averaging. We prove this result under two general settings--with and without spatial embeddings--and only mild technical restrictions on the form of microscopic nonlinearity, noise distributions, and network connectivity patterns. Notably, however, the limit dynamics are in general linear time-varying (LTV).
\textit{Our third contribution} then consists of extending the above results to scenarios where, under additional stationarity-related assumptions on the microscopic dynamics, the spatially-averaged limit dynamics are linear and time-invariant (LTI).
Finally, \textit{our fourth contribution} pertains to the finite-sample analysis of the rate at which averaged dynamics converge to linearity. In essence, our analysis shows that in the extreme case of independent subsystems (where averaging is most effective), convergence occurs at the parametric rate ($O(1/\sqrt N)$), and the convergence slows down as the correlations between subsystems strengthen (with the opposite extreme of no convergence if all the subsystems are perfectly correlated).
Notably, among these technical contributions the only one that existed in our preliminary version~\cite{ahmed2022linearizing} is the third. The first and the fourth sets of contributions are completely novel, while our second contribution is now proven under significantly more general conditions--particularly, interconnected and heterogeneous systems--whereas its preliminary version was shown only for homogeneous and disconnected (i.i.d.) subsystems.
To our knowledge, our results provide the first theoretical account of emergent linearity in complex dynamical systems, and open the door to a broad range of subsequent empirical and theoretical studies.

\section{Notation}

We use $\real$ to denote the set of real numbers.
For a matrix $\Abf \in \real^{m \times n}$, $\Abf^\dagger \in \real^{n \times m}$ denotes its pseudoinverse. Throughout this work, all probabilities are defined on measurable spaces consisting of a Euclidean space (or a subset thereof) and the associated Borel $\sigma$-algebra. Hence, when clear from the context, the space over which each probability is defined is omitted. $\E[\cdot]$ and $\P\{\cdot\}$ denote expectation and probability, respectively. For two random variables $\xi$ and $\eta$, $\var(\xi)$ and $\rho(\xi, \eta)$ denote the variance of $\xi$ and the Pearson correlation coefficient between $\xi$ and $\eta$, respectively. For two random vectors $\xib \in \real^n$ and $\etab \in \real^m$, $\cov(\xib, \etab) \in \real^{n \times m}$ is the covariance between them, and $\cov(\xib) = \cov(\xib, \xib)$. For sequences of random variables, we use both $\stackrel{d}{\to}$ and $\dlim$ to denote their convergence in distribution and $\stackrel{a.s.}{\to}$ to denote almost sure convergence.

\section{Preliminaries: Mixing Sequences}\label{sec:prelims}

In this section we review some fundamental notions and properties of mixing sequences which form the theoretical basis for the ensuing discussion.
In essence, for a discrete-time stochastic process $\xi_i, i = 1, 2, \dots$ different notions of `mixing' characterize the case when the statistical dependence between $\xi_i$ and $\xi_j$ diminishes as $|i - j|$ increases. As such, mixing conditions generalize the notion of a pairwise independent (a.k.a. white) sequence to one in which nearby elements can be dependent but their dependence decays as the distance between them grows. Various versions of mixing sequences have been proposed, corresponding to different measures of dependence which has to decay with distance. One of the most practical and empirically verifiable versions is that of $\rho$-mixing, as defined next.

\begin{definition}\longthmtitle{$\rho$-mixing sequence~\cite{bradley2005basic}}\label{def:rho-mixing}
Consider a sequence of random variables $\xi_1, \xi_2, \dots$ in a probability space $(\Omega, \mathcal{F}, \mathcal{P})$
and define
\begin{align}\label{eq:rho_n}
    \varrho(n) = \sup_{i} \ \rho(\sigma(\xi_1, \dots, \xi_i), \sigma(\xi_{i+n}, \xi_{i+n+1}, \dots)),
\end{align}
where $\sigma(\cdot)$ denotes the smallest $\sigma$-algebra of $\Omega$ generated by a set of random variables, and, for any two $\sigma$-algebras $\mathcal{A}$ and $\mathcal{B}$,
\begin{align}\label{eq:rhoAB}
    \rho(\mathcal{A}, \mathcal{B}) = \sup \frac{|\cov(y, z)|}{\var(y)^{1/2} \var(z)^{1/2}}
\end{align}
where the supremum is taken over all pairs of squared-integrable random variables $y$ and $z$ such that $y$ is $\mathcal{A}$-measurable and $z$ is $\mathcal{B}$-measurable.
The sequence $\xi_1, \xi_2, \dots$ is then $\rho$-mixing if
\begin{align*}
    \varrho(n) \to 0 \quad \text{as} \quad n \to \infty. \eqoprocend
\end{align*}
\end{definition}

Clearly, any i.i.d. sequence is $\rho$-mixing and the latter is a generalization of the former. Thus, many properties of i.i.d. sequences, such as the laws of large numbers~\cite{fazekas2001general, andrews1988laws} and the central limit theorem (CLT)~\cite{ekstrom2014general, functional_clt}, have been generalized to mixing sequences as long as the decay of dependence is sufficiently fast. Of particular relevance to this work is the CLT for $\rho$-mixing sequences, as presented below. But we first need a technical definition, as follows.

\begin{definition}\longthmtitle{Slowly-varying function}\label{def:slowly-varying-function}
A function $h: \intpos \to \real$ is called \emph{slowly varying} if
\begin{align*}
    \lim_{N \to \infty} \frac{h(kN)}{h(N)} = 1
\end{align*}
for all $k \in \intpos$. \oprocend
\end{definition}

Intuitively, if $h(N)$ is slowly varying it implies that the asymptotic behavior of $h(N)$ becomes insensitive to constant multiplicative factors, resulting in a very slow rate of growth/decay.
Some examples of slowly-varying functions are constant functions $h(N) = c$, logarithm functions such as $h(N) = \log N$ or $h(N) = \log \log N$, and bounded functions with nonzero constant limits such as $h(x) = 1 + \frac{\sin(x)}{x}$.

We are now ready to present the following result, which is a generalization of the standard CLT for i.i.d. sequences~\cite[Thm 3.4.1]{durrett2019probability} to $\rho$-mixing sequences of random variables.

\begin{proposition}\longthmtitle{CLT for $\rho$-mixing sequences~\cite[Thm B]{functional_clt}}\label{prop:clt}
Consider a $\rho$-mixing sequence of random variables $\xi_1, \xi_2, ...$. %
and define its cumulative sum and cumulative variance as
\begin{align}\label{eq:sum_var}
    S_N = \sum_{i=1}^{N} \xi_i, \quad \sigma_N^2 = \var(S_N).
\end{align}
Assume
\begin{subequations}
\begin{align}
\label{prop:cond_1} &\E[\xi_i] = 0 \ \text{and} \ \var(\xi_i) < \infty \ \text{for all i}, \qquad\qquad\qquad\quad
\\
\label{prop:cond_3} &\sup_{M \ge 0, N \ge 1} \ \frac{1}{\sigma_N^2} \ \E[(S_{M+N} - S_M)^2] < \infty,
\\
\label{prop:cond_2} &h(N) = \frac{\sigma_{N}^2}{N} \ \text{is a slowly-varying function.}
\end{align}
\end{subequations}
Then
\begin{align}\label{eq:rho-mixing-clt}
    \frac{S_N}{\sigma_N} \stackrel{d}{\to} \Nc(0, 1) \quad \text{as} \quad N \to \infty. \quad \oprocendsymbol
\end{align}
\end{proposition}

It is instructive to compare Proposition~\ref{prop:clt} with the standard CLT for i.i.d. sequences~\cite[Thm 3.4.1]{durrett2019probability}. Both versions assume that all $\xi_i$ have finite mean and variance; they slightly differ in that Proposition~\ref{prop:clt} assumes $\xi_i$'s are already mean-subtracted while the standard version mean-subtracts the sum.
Condition~\eqref{prop:cond_3} asks that the ratio of the variance of shifted partial sums to that of initial partial sums remains bounded. This is a generalization of the `identically distributed' assumption in the standard version, in which case the left hand side of~\eqref{prop:cond_3} is always 1.
Finally, $\rho$-mixingness is a generalization of independence, as noted earlier. Importantly, however, condition~\eqref{prop:cond_2} ensures that the sequence $\{\xi_i\}$ is not any $\rho$-mixing sequence, but one in which the correlations between neighboring elements decays sufficiently fast. This is a subtle point and we will get back to it later in Section~\ref{sec:s_avg_1d} when we define residual factors. For now, note
the normalization by $\sigma_N$ in~\eqref{eq:rho-mixing-clt} vs. the standard normalization by $\sqrt N$ for i.i.d. sequences. This difference stems from the fact that the growth rate of the cumulative variance ($\sigma_N^2$) of $\rho$-mixing sequences can vary, whereas $\sigma_N^2$ always grows as $N$ for i.i.d. sequences.

In what follows, we will extend Proposition~\ref{prop:clt} in a number of ways, and use the results to prove the Gaussianity of joint state-noise distributions under spatial averaging, using which we will prove the linearity of spatially-averaged dynamics.

\section{Linearizing Effect of Spatial Averaging on Sequences of Dynamical Systems}\label{sec:s_avg_1d}

\subsection{Problem Formulation}

Consider a heterogeneous population of $N$ dynamical subsystems, %
where each has the general discrete-time nonlinear form
\begin{align}\label{eq:ode}
     \xbf_i(t+1) &= f_i(\xbf_{\Nc_i}(t), \wbf_i(t)), \qquad i = 1, 2, \dots N,
     \\
     \notag \xbf(0) &= \begin{bmatrix}
     \xbf_1(0)^T & \cdots & \xbf_N(0)^T
     \end{bmatrix}^T \sim p_0,
     \\
     \notag \wbf(t) &= \begin{bmatrix}
     \wbf_1(t)^T & \cdots & \wbf_N(t)^T
     \end{bmatrix}^T \sim p_w(t).
\end{align}
$\xbf_i(t) \in \real^n$ is the state of subsystem $i$ with initial joint distribution
$p_0$, and $\wbf_i(t) \in \real^m$ is the noise process of subsystem $i$ with the joint distribution
$p_w(t)$. The vector
\begin{align}\label{eq:x_Ni}
    \xbf_{\Nc_i}(t) =
    \begin{bmatrix}
        \xbf_{\max\{1,i-\tau\}}^T & \!\cdots\! & \xbf_i^T & \!\cdots\! & \xbf_{\min\{i+\tau,N\}}^T
    \end{bmatrix}^T
\end{align}
denotes the state of all subsystems on which $f_i$ \textit{can} depend (i.e., a superset of the in-neighbor set of node $i$). $\tau < \infty$ can be arbitrarily large, but finite, and controls the spatial range of neighborhoods, such that subsystems $i$ and $j$ can (but do not need to) be neighbors if $|i - j| \le \tau$.

The problem we tackle in this section, motivated by our prior empirical observations~\cite{nozari2023macroscopic} and those of others~\cite{sani2021modeling,schulz2019deep,yang2019dynamic}, is as follows.

\begin{problem}\longthmtitle{Linearizing Effect of Spatial Averaging on Sequences of Dynamical Systems}\label{prob1}
    Consider a heterogeneous population of nonlinear dynamical systems described by~\eqref{eq:ode}, and define the population's average state vector as
    \begin{align}\label{eq:s_avg_state}
        \bar \xbf(t) = \frac{1}{\phi(N)} \sum_{i=1}^{N} \xbf_i(t) - \mathbb{E}[\xbf_i(t)],
    \end{align}
    where $\phi(N)$ is a normalization factor. Prove, under appropriate assumptions and choice of $\phi(N)$, that
    the dynamics of $\bar \xbf(t)$ becomes asymptotically linear as $N \to \infty$. \oprocend
\end{problem}

The normalization factor $\phi(N)$ plays the same role that $\sigma_N$ plays in~\eqref{eq:rho-mixing-clt}, but it is simpler as it only captures the \textit{growth rate} of $\sigma_N$ and not $\sigma_N$ itself. For comparison, in the standard CLT $\phi(N) = \sqrt N$ for all i.i.d. sequences while $\sigma_N$ varies from one sequence to another.  %
As we move to multivariate processes (see below) this distinction will become bolder. A naive extension of~\eqref{eq:rho-mixing-clt} would require computing the exact covariance matrix of $S_N$, which often lacks a closed-form solution, as well as normalizing $S_N$ by $\cov(S_N)^{-\frac{1}{2}}$, which would inter-mix different dimensions and lose the notion of an `average'.

In the rest of this section we will first lay the statistical foundation of our framework in Section~\ref{sub_sec:rho-mixing-seq}, where we extend and further characterize the properties of $\rho$-mixing sequences. In Section~\ref{sub_sec:s_avg_1d} we then use this foundation to present our main result of this section, i.e., the %
solution to Problem~\ref{prob1}. Under additional stationarity assumptions, we further strengthen this result to prove linear \textit{time-invariant} (LTI) limit dynamics, and finally characterize the rate of this convergence to linearity in Section~\ref{sub_sec:finite_err} for the special case of i.i.d. sequences.

\subsection{Multivariate $\rho$-Mixing Sequences}\label{sub_sec:rho-mixing-seq}

In this subsection we extend Definition~\ref{def:rho-mixing} to vector-valued stochastic processes, and then prove some of their properties, including a more practical CLT.
These results will play a central role in our proof of the linearizing effect of spatial averaging in Theorem~\ref{thm:s_avg}.

The following is a natural generalization of Definition~\ref{def:rho-mixing}.

\begin{definition}\longthmtitle{Multivariate $\rho$-mixing sequences}\label{def:mv-rho-mixing}
Consider a sequence of random vectors $\xib_1, \xib_2, \dots \in \real^q$ in a probability space $(\Omega, \mathcal{F}, \mathcal{P})$,
and define
\begin{align}\label{eq:rho_n_mv}
    \varrho(n) = \sup_{j} \ \rho(\sigma(\xib_1, \dots, \xib_j), \sigma(\xib_{j+n}, \dots)),
\end{align}
where$\rho(\Ac, \Bc)$ for two $\sigma$-algebras $\Ac$ and $\Bc$ is defined in~\eqref{eq:rhoAB}.
The sequence $\xib_1, \xib_2, \dots$ is called $\rho$-mixing if
\begin{align}\label{eq:rhoto0}
    \varrho(n) \to 0 \quad \text{as} \quad n \to \infty.
\end{align}
\end{definition}

Before we can present the generalized CLT for multivariate $\rho$-mixing sequences we also need the following definition. This will help, in particular, in determining the appropriate normalization factor in Problem~\ref{prob1}.

\begin{definition}\longthmtitle{Residual factor}\label{def:resfact}
Let $\xib_1, \xib_2, \dots$ be a $\rho$-mixing sequence and $S_N = \sum_{i = 1}^N \xib_i$. The function $h(N)$ is called a \emph{residual factor} for $\{\xib_i\}_{i = 1}^\infty$ if
\begin{align}\label{eq:resfact}
\lim_{N \to \infty} \frac{1}{N h(N)} \cov(S_N) < \infty,
\end{align}
i.e., the limit exists and is finite. \oprocend
\end{definition}

To put this definition in perspective, if $\{\xib_i\}_{i = 1}^\infty$ is i.i.d. then $\cov(S_N) = N \cov(\xib_1)$ and so $h(N) = 1$ is a residual factor. Note, also, that residual factors are not unique. In the i.i.d. case, $h(N) = c$ is also a residual factor for any constant $c$. Further, note that Definition~\ref{def:resfact} does not require the limit %
to be nonzero. Thus, $h(N) = N$, $\log(N)$, and $e^N$ are all valid residual factors for an i.i.d. sequence. Finally, and related to the last point, a residual factor always exists for any $\rho$-mixing sequence. The following Theorem formalizes this existence.

\begin{theorem}\longthmtitle{Existence of residual factors}\label{thm:resfact-exist}
Let $\xib_1, \xib_2, \dots \in \real^q$ be a $\rho$-mixing sequence and assume that $\var((\xib_i)_\ell) \le \bar \sigma^2$ for all $i, \ell$ and some $\bar \sigma < \infty$. Then, $h(N) = N$ is a residual factor for $\{\xib_i\}_{i = 1}^\infty$.
\end{theorem}
\begin{proof}
Let $S_N = \sum_{i = 1}^N \xib_i$ and let $\varrho(n)$ be as in~\eqref{eq:rho_n_mv}. Then, for any $\ell = 1, \dots, q$,
\begin{align}\label{eq:varub1}
\notag \var((S_N)_\ell) &= \sum_{i=1}^N \var((\xib_i)_\ell) + 2 \sum_{1 \leq i < j \leq N} \cov((\xib_i)_\ell, (\xib_j)_\ell)
\\
&\le N \bar \sigma^2 + 2 \bar \sigma^2 \sum_{1 \leq i < j \leq N} \varrho(|i-j|).
\end{align}
Changing the order of summation, we have:
\begin{align}\label{eq:varub2}
\notag \sum_{1 \leq i < j \leq N} \varrho(|i-j|) &= \sum_{k=1}^{N-1} \sum_{i=1}^{N-k} \varrho(k) = \sum_{k=1}^{N-1} (N-k)\varrho(k)
\\
&\le N \sum_{k=1}^N \varrho(k).
\end{align}
Assume, without loss of generality, that $\varrho(x)$ is defined and infinitely differentiable for all real-valued $x \in \realpos$. This is without loss of generality since $\varrho(n)$ can always be replaced by an upper bound that satisfies this assumption and still decays to 0 as $n \to \infty$. Then, using the Euler–Maclaurin formula~\cite{graham1994concrete},
\begin{align}\label{eq:varub3}
\sum_{k=1}^N \varrho(k) = O\Big(\int_1^N \varrho(x) dx\Big)
\end{align}
Combining~\eqref{eq:varub1}-\eqref{eq:varub3}, we get
\begin{align}\label{eq:twoterms}
\lim_{N \to \infty} \frac{1}{N^2} \var((S_N)_\ell) \le \lim_{N \to \infty} \frac{N \bar \sigma^2}{N^2} + \frac{2\bar \sigma^2 C}{N} \int_1^N \varrho(x) dx
\end{align}
where $C$ is a bounding constant from~\eqref{eq:varub3}. Limit of the first term is clearly 0. The second limit has two possibilities. If the non-negative quantity
\begin{align}\label{eq:limint}
\int_1^\infty \varrho(x) dx
\end{align}
is finite ($\varrho(n)$ decays fast), then the second limit on the right hand side of~\eqref{eq:twoterms} is also clearly 0. If, on the other hand, the limit in~\eqref{eq:limint} is infinite ($\varrho(n)$ decays slowly), then by the L'Hopital's rule and~\eqref{eq:rhoto0}, still
\begin{align*}
\lim_{N \to \infty} \frac{\int_1^N \varrho(x) dx}{N} = \lim_{N \to \infty} \frac{\varrho(N)}{1} = 0.
\end{align*}
Therefore, either way, the limits in~\eqref{eq:twoterms} are all 0.
Since this is true for all $\ell = 1, \dots, q$, we also have
\begin{align*}
\lim_{N \to \infty} &\frac{1}{N^2} |\cov((S_N)_\ell, (S_N)_m)|
\\
&\le \lim_{N \to \infty} \frac{\sqrt{\var((S_N)_\ell)}}{N} \frac{\sqrt{\var((S_N)_m)}}{N} = 0,
\end{align*}
for all $\ell, m = 1, \dots, q$. Put together, $\lim_{N \to \infty} \frac{1}{N^2} \cov(S_N) = 0$, completing the proof.
\end{proof}

A corollary to Theorem~\ref{thm:resfact-exist} is that any $h(N) \ge N$ is also a valid residual factor for any $\rho$-mixing sequence with uniformly-bounded variance. However, residual factors are practically useful only when they grow at the same rate as $\cov(S_N)$ itself and the limit in~\eqref{eq:resfact} is nonzero. This ``smallest" residual factor is what we will use in the sequel, even though the forthcoming results will still be correct, but not necessarily useful, for all valid residual factors.

Going back to Proposition~\ref{prop:clt}, note also the similarity between~\eqref{prop:cond_2} and~\eqref{eq:resfact}, where $h(N)$ in~\eqref{prop:cond_2} is clearly a residual factor for $\{\xi_k\}_{i = 1}^\infty$ in Proposition~\ref{prop:clt}. As we will see later, whether a residual factor of a $\rho$-mixing sequence is slowly-varying remains to be a critical condition for validity of CLT (and in turn asymptotic linearity) for that sequence.

To better clarify the notion of a residual factor and how it depends on the statistics of $\{\xib_i\}_{i = 1}^\infty$, in the following result we provide closed-form expressions for ``smallest" residual factors of a few representative $\rho$-mixing sequences.

\begin{theorem}\longthmtitle{Relationship between decay r of correlations and growth of residual factors}\label{thm:hN}
Let $\xi_1, \xi_2, \dots \in \real$ be a $\rho$-mixing sequence of random variables, and let $\rho(n)$ be as in~\eqref{eq:rho_n_mv}. Assume, for simplicity, that $\var(\xi_i) = \sigma^2$ for all $i$ and $\rho(\xi_i, \xi_j) = \varrho(|i-j|)$. Then the following relationships exist, where $h(N)$ is a residual factor in each case:
\begin{enumerate}
\item if $\varrho(n) = \frac{1}{n}$ then $h(N) = \log(N)$;
\item if $\varrho(n) = \frac{1}{\log(1+n)}$ then $h(N) = \frac{N}{\log(N)}$;
\item if $\varrho(n) = \frac{1}{n^p}, p \in (0, 1)$, then $h(N) = N^{1-p}$.
\end{enumerate}
\end{theorem}
\begin{proof}
In all cases, let $S_N = \sum_{i=1}^N \xi_i$.
\begin{enumerate}[wide]
\item Similar to~\eqref{eq:varub1} and~\eqref{eq:varub2}, we have
\begin{align*}
\text{Var}(S_N) &= N \sigma^2 + 2 \sigma^2 \sum_{1 \leq i < j \leq N} \frac{1}{|i-j|},
\\
\sum_{1 \leq i < j \leq N} \frac{1}{|i-j|} &= %
\sum_{k=1}^{N-1} \frac{N-k}{k} = N \sum_{k=1}^{N-1} \frac{1}{k} - (N-1).
\end{align*}
For large $N$, the harmonic sum $\sum_{k=1}^{N-1} \frac{1}{k}$ is approximately $\log(N) + \gamma$, where $\gamma$ is the Euler's constant. Thus, for large $N$,
\begin{align*}
\sum_{1 \leq i < j \leq N} \frac{1}{|i-j|} &\simeq N (\log(N) + \gamma) - (N-1)
\\
&= N \log(N) + \text{lower order terms},
\end{align*}
and
\begin{align*}
\text{Var}(S_N) &\simeq N \sigma^2 + 2 \sigma^2 N \log(N) + \text{lower order terms} \\
&= 2 \sigma^2 N \log(N) + \text{lower order terms}.
\end{align*}
Therefore, for large $N$,
\begin{align*}
\frac{\var(S_N)}{N} \simeq 2 \sigma^2 \log(N) + \text{lower order terms}
\end{align*}
and so $h(N) = \log(N)$ is a residual factor. This is also a slowly-varying function, since for $k \in \intpos$,
\begin{align*}
\lim_{N \to \infty} \frac{h(k N)}{h(N)} = \lim_{N \to \infty} \frac{ 2 \sigma^2 \log(k N)}{2 \sigma^2 \log(N)} = 1.
\end{align*}
In other words, the correlation decay rate of $\frac{1}{|i - j|}$ is ``sufficiently fast" so that $\var(S_N)$ grows sufficiently slowly and $h(N)$ becomes slowly varying. As we will see in Theorem~\ref{thm:multivar_clt_rho_mixing}, this allows for CLT to hold for this sequence.

\item Similar to case (i),
\begin{align*}
\text{Var}(S_N) %
&= N \sigma^2 + 2 \sigma^2 \sum_{1 \leq i < j \leq N} \frac{1}{\log(1+|i-j|)}
\end{align*}
and
\begin{align}\label{eq:expression2}
\notag &\sum_{1 \leq i < j \leq N} \frac{1}{\log(1+|i-j|)} %
= \sum_{k=1}^{N-1} \frac{N-k}{\log(1+k)} = \sum_{k=2}^N \frac{N\!+\!1\!-\!k}{\log(k)}
\\
&\qquad= (N+1)\sum_{k = 2}^N \frac{1}{\log(k)} - \sum_{k = 2}^N \frac{k}{\log(k)}.
\end{align}
For large $N$, this can be approximated (with the same growth rate) using the Euler-Maclaurin formula~\cite{graham1994concrete} as
\begin{align}\label{eq:expression1}
(N+1)\int_2^N \frac{1}{\log(x)}dx - \int_2^N \frac{x}{\log(x)} dx.
\end{align}
By definition, $\int_2^N \frac{1}{\log(x)}dx = \li(N) - \li(2)$ where $\li(u) = \int_0^u \frac{dx}{\log(x)}$ is the logarithmic integral function. To compute the second integral, we can use the change of variables $u = \log(x)$ and $v = 2u$ so that
\begin{align*}
\int_2^N \frac{x}{\log(x)} dx &= \int_{\log(2)}^{\log(N)} \frac{e^u}{u} e^u du = \int_{2\log(2)}^{2\log(N)} \frac{e^v}{v} dv
\\
&= \Ei(\log(N^2)) - \Ei(\log(4))
\end{align*}
where $\Ei(u) = \int_{-\infty}^u \frac{e^x}{x} dx$ is the exponential integral function. Since $\li(u) = \Ei(\log(u))$, \eqref{eq:expression1} further simplifies to
\begin{align}\label{eq:expression3}
N \li(N) - \li(N^2) + \li(N) - (N+1) \li(2) + \li(4)
\end{align}
For large $N$,
\begin{align*}
\li(N) = \frac{N}{\log(N)} \Big[1 + O\big(\frac{1}{\log(N)}\big)\Big] \simeq \frac{N}{\log(N)}
\end{align*}
Therefore, \eqref{eq:expression3} is approximately equal to
\begin{align*}
&\frac{N^2}{\log(N)} - \frac{N^2}{\log(N^2)} + \frac{N}{\log(N)} - (N+1) \li(2) + \li(4)
\\
&= \frac{N^2}{2\log(N)} + \text{lower order terms}.
\end{align*}
Combining these results, for large $N$ we get
\begin{align*}
\frac{\text{Var}(S_N)}{N} &\simeq \sigma^2 + 2 \sigma^2 \frac{N}{2\log(N)} + \text{lower order terms}
\\
&\simeq \frac{\sigma^2 N}{\log(N)} + \text{lower order terms}.
\end{align*}
and so $h(N) = \frac{N}{\log(N)}$ is a residual factor.
This is not a slowly-varying function, however, since for $k \in \intpos$,
\begin{align*}
\lim_{N \to \infty} \frac{h(k N)}{h(N)} = \lim_{N \to \infty} \frac{\frac{\sigma^2 k N}{\log(k N)}}{\frac{\sigma^2 N}{\log(N)}} = k.
\end{align*}
Therefore in this case, unlike case (i), correlations decay too slowly. As a result, $\var(S_N)$ grow too fast and $h(N)$ is not slowly-varying. As we will see in Theorem~\ref{thm:multivar_clt_rho_mixing}, this prevents CLT to hold for this sequence, despite being $\rho$-mixing.

\item This case falls in between cases (i) and (ii) in terms of the decay rate of $\rho(n)$.
Proceeding as before,
\begin{align*}
\text{Var}(S_N) = N \sigma^2 + 2 \sigma^2 \sum_{1 \leq i < j \leq N} \frac{1}{|i-j|^p}
\end{align*}
and
\begin{align*}
\sum_{1 \leq i < j \leq N} \frac{1}{|i-j|^p} = \sum_{k=1}^{N-1} \frac{N-k}{k^p} = N\sum_{k=1}^{N-1} k^{-p} - \sum_{k=1}^{N-1} k^{1-p}
\end{align*}
Using the Euler-Maclaurin formula, we get for large $N$,
\begin{align*}
\sum_{k=1}^{N-1} k^{-p} &\simeq \int_1^{N-1} x^{-p} dx = \frac{(N-1)^{1-p} - 1}{1-p}
\\
\sum_{k=1}^{N-1} k^{1-p} &\simeq \int_1^{N-1} x^{1-p} dx = \frac{(N-1)^{2-p} - 1}{2-p}
\end{align*}
Thus, for large $N$,
\begin{align*}
\sum_{1 \leq i < j \leq N} \frac{1}{|i-j|^p} &\simeq \frac{N(N-1)^{1-p}}{1-p} - \frac{(N-1)^{2-p}}{2-p}
\\
&= \frac{(N-1)^{2-p}}{(1-p)(2-p)} + \frac{(N-1)^{1-p}}{1-p}
\end{align*}
and so
\begin{align*}
\frac{\text{Var}(S_N)}{N} &\simeq \sigma^2 + \frac{2 \sigma^2}{N} \Big[\frac{(N-1)^{2-p}}{(1-p)(2-p)} + \frac{(N-1)^{1-p}}{1-p}\Big]
\\
&=  \frac{2 \sigma^2}{(1-p)(2-p)} \frac{(N-1)^{2-p}}{N} + \text{lower order terms}.
\end{align*}
Therefore, $h(N) = \frac{N^{2-p}}{N} = N^{1-p}$ is a residual factor.
This $h(N)$ is not a slowly-varying function, since for $k \in \intpos$,
\begin{align*}
\lim_{N \to \infty} \frac{h(k N)}{h(N)} = \lim_{N \to \infty} \frac{(kN)^{1-p}}{N^{1-p}} = k^{1-p}.
\end{align*}
In other words, even the polynomial decay rate of $\frac{1}{|i-j|^p}, p \in (0, 1)$ is too slow for $h(N)$ to be slowly varying and for CLT to hold.
\end{enumerate}
\vspace*{-13pt}
\end{proof}

In summary, it follows from Theorem~\ref{thm:hN} that
$\rho(n) \propto \frac{1}{n}$
is the rough ``boundary" between decay rates that are sufficiently fast %
and those that are too slow. %

Before we are ready to present the generalized CLT for multivariate $\rho$-mixing sequences, we need to prove another property of $\rho$-mixing sequences, namely, that sequences formed through transformations of $\rho$-mixing sequences are also $\rho$-mixing. %

\begin{lemma}\longthmtitle{Sequence formed through transformation of $\rho$-mixing sequence is $\rho$-mixing}\label{lem:rho_mixing_seq}
Consider a sequence of random vectors $\xib_1, \xib_2, \dots \in \real^q$ in a probability space $(\Omega, \mathcal{F}, \mathcal{P})$. Let $\tau, k \ge 0$ be arbitrary fixed integers and assume that for any $i \ge 1$, $h_i$ is a measurable function from $\real^{q\min\{\tau+i, 2\tau+1\}}$ to $\real^k$. Define
\begin{align*}
    \zetab_i = h_i(\xib_{\min\{1,i-\tau\}}, \dots, \xib_i, \dots, \xib_{i+\tau}), \qquad i \ge 1.
\end{align*}
If the sequence $\xib_1, \xib_2, \dots$
is $\rho$-mixing, then so is $\zetab_1, \zetab_2, \dots$.
\end{lemma}
\begin{proof}
Without loss of generality let $n \ge \tau + 1$. Then, by the definition of the $\sigma$-algebra generated by a collection of random variables and the measurability of all $h_i$, we get
\begin{align*}
\sigma(\zetab_1, \dots, \zetab_j) \subseteq \sigma(\xib_1, \dots, \xib_{j+\tau}),
\end{align*}
and
\begin{align*}
\sigma(\zetab_{j+n}, \dots) \subseteq \sigma(\xib_{j+n-\tau}, \dots),
\end{align*}
for all $j \ge 1$. Therefore, by~\eqref{eq:rhoAB},
\begin{align}\label{eq:sandwich}
0 \le \rho_\zeta(n) \le \rho_\xi(n - 2\tau), \qquad n \ge \tau+1.
\end{align}
where $\rho_\zeta(n)$ and $\rho_\xi(n)$ are defined as in~\eqref{eq:rho_n} for the respective sequences $\{\zetab_i\}$ and $\{\xib_i\}$. The theorem then follows by letting $n \to \infty$ in~\eqref{eq:sandwich} and using the mixing property of $\{\xib_i\}$.
\end{proof}

The next result extends and improves the CLT in Proposition~\ref{prop:clt} to the case of multivariate $\rho$-mixing sequences. %

\begin{theorem}\longthmtitle{Multivariate CLT for $\rho$-mixing sequence}\label{thm:multivar_clt_rho_mixing}
Let $\xib_1, \xib_2, ... \in \real^q$ be a $\rho$-mixing sequence where each $\xib_i$ has zero mean and finite variance. Let $h(N)$ be a residual factor for $\{\xib_i\}_{i = 1}^\infty$ and assume that it is slowly-varying.
Also, assume that for all $\thetab \neq 0$,
\begin{align}\label{eq:st-var-growth}
\sup_{M, N} \frac{\thetab^T \cov\Big(\sum_{i = M+1}^{M+N} \xib_i \Big) \thetab}{\thetab^T \cov\Big(\sum_{i = 1}^N \xib_i \Big) \thetab} < \infty.
\end{align}
Then, the ``mean" variable
\begin{align}\label{eq:xibar}
    \bar \xib = \frac{1}{\sqrt{Nh(N)}} \sum_{i=1}^N \xib_i,
\end{align}
satisfies
\begin{align}\label{eq:rho-mixing-clt-fin-var}
    \bar \xib \stackrel{d}{\longrightarrow} \Nc(\zeros, \Sigmab_{\bar \xi}^*) \quad \text{as} \quad N \to \infty,
\end{align}
where (cf.~\eqref{eq:resfact})
\begin{align}\label{eq:finite_cov}
     \Sigmab_{\bar \xi}^* = \lim_{N \to \infty} \frac{1}{N h(N)} \cov\Big(\sum_{i=1}^N \xib_i\Big).
\end{align}
\end{theorem}
\begin{proof}
Let $\thetab \neq \zeros$ be an arbitrary constant vector. By Lemma~\ref{lem:rho_mixing_seq}, the sequence of random variables $\{\thetab^T \xib_i\}_{i \ge 1}$ is $\rho$-mixing. Define the cumulative sum and cumulative variance of this sequence as
\begin{align}\label{eq:new_sum_var}
  S_N = \sum_{i=1}^N \thetab^T \xib_i, \qquad \sigma_N^2 = \var(S_N).
\end{align}
The sequence $\{\thetab^T \xib_i\}_{i \ge 1}$ follows~\eqref{prop:cond_1} by assumption. It also follows~\eqref{prop:cond_3} due to~\eqref{eq:st-var-growth}. We need to show that it also satisfies~\eqref{prop:cond_2} before we can use the CLT in Proposition~\ref{prop:clt}.

Note that
\begin{align*}
 \sigma_N^2 = \thetab^T \cov\Big(\sum_{i=1}^N \xib_i\Big) \thetab
\end{align*}
so we get
\begin{align}\label{eq:new_cond_2}
\lim_{N \to \infty} \frac{\sigma_N^2}{Nh(N)} &= \thetab^T \Sigmab_{\bar \xi}^* \thetab.
\end{align}
First consider the case where $\thetab^T \Sigmab_{\bar \xi}^* \thetab \neq 0$.
Using the fact that $h(N)$ is slowly varying, we have, for any $k \in \intpos$,
\begin{align*}
\lim_{N \to \infty} \frac{\frac{\sigma_{kN}^2}{kN h(kN)} h(kN)}{\frac{\sigma_N^2}{N h(N)} h(N)}
&= \frac{ \lim\limits_{N \to \infty} \frac{\sigma_{kN}^2}{kN h(kN)}}{ \lim\limits_{N \to \infty} \frac{\sigma_N^2}{N h(N)}} \cdot \lim_{N \to \infty} \frac{h(kN)}{h(N)}
\\
&= \frac{\thetab^T \Sigmab_{\bar \xi}^* \thetab}{\thetab^T \Sigmab_{\bar \xi}^* \thetab} \cdot 1 = 1.
\end{align*}
This proves that
\begin{align}\label{eq:new_slowly_var}
\tilde h(N) = \frac{\sigma_N^2}{N}
\end{align}
is slowly varying, and thus $\sigma_N^2$ satisfies~\eqref{prop:cond_2}.
Therefore, Proposition~\ref{prop:clt} holds and, from~\eqref{eq:rho-mixing-clt} and~\eqref{eq:new_cond_2}, as $N \to \infty$
\begin{align}\label{eq:clt-result}
\notag \frac{S_N}{\sqrt{Nh(N)}} = \frac{S_N}{\sigma_N} \cdot \frac{\sigma_N}{\sqrt{Nh(N)}} &\stackrel{d}{\to} \Nc(0,1) \cdot \sqrt{\thetab^T \Sigmab_{\bar \xi}^* \thetab} \\
&= \Nc(0, \thetab^T \Sigmab_{\bar \xi}^* \thetab)
\end{align}
This means that
\begin{align}\label{eq:cramer-wold-device-eq-mult}
\thetab^T\bar \xib \stackrel{d}{\to} \thetab^T\bar \xib^*, \quad \text{as} \quad N \to \infty,%
\end{align}
where $\bar \xib^*$ is a multivariate random variable distributed as $\Nc\left(\zeros,  \Sigmab_{\bar \xib}^*\right)$.

On the other hand, if $\thetab^T \Sigmab_{\bar \xi}^* \thetab = 0$ (including when $\thetab = \zeros$), then \eqref{eq:cramer-wold-device-eq-mult} follows immediately without the need to use any CLT since $\thetab^T\bar \xib^* = 0$ and
\begin{align*}
\lim_{N \to \infty} \var(\thetab^T \bar \xib) = \lim_{N \to \infty} \frac{\sigma_N^2}{Nh(N)} \stackrel{\text{\eqref{eq:new_cond_2}}}{=} \thetab^T \Sigmab_{\bar \xi}^* \thetab = 0.
\end{align*}
Therefore, \eqref{eq:cramer-wold-device-eq-mult} holds for all $\thetab \in \real^q$. According to Cramer-Wold device theorem~\cite[Thm 3.9.5]{durrett2019probability}, it then follows that
\begin{align}\label{eq:zbar_stozstar_s}
\bar \xib \stackrel{d}{\longrightarrow} \bar \xib^*  \quad \text{as} \quad  N \to \infty,
\end{align}
completing the proof.
\end{proof}

Theorem~\ref{thm:multivar_clt_rho_mixing} provides the foundation for the asymptotic analysis of spatially averaged dynamics formulated in Problem~\ref{prob1}, which we undertake next.

 \subsection{Asymptotic Linearity Under Spatial Averaging}\label{sub_sec:s_avg_1d}

In this section we use the results in Section~\ref{sub_sec:rho-mixing-seq} to tackle Problem~\ref{prob1}.
Throughout this section, we make the following assumptions on the population dynamics~\eqref{eq:ode}. For ease of notation, for all $i, t$ let
\begin{align}\label{eq:rho-mixing-1d}
   \ybf_i(t) = \begin{bmatrix}
        \xbf_i(t)^T & \wbf_i(t)^T
    \end{bmatrix}^T,
\end{align}
and
\begin{align}\label{eq: conc_states}
    \zbf_i(t) =  \begin{bmatrix}
    \xbf_i(t+1)^T & \xbf_i(t)^T & \wbf_i(t)^T
    \end{bmatrix}^T.
\end{align}

\begin{assumption}\longthmtitle{Standing assumptions}\label{assum}
    \begin{enumerate}[label=\rm(A\arabic*), leftmargin=25pt]
        \item\label{assum:w} The noise processes are
        zero mean %
        and finite-variance, i.e., there exists $C < \infty$ such that for all $i$ and $t$,
        \begin{align*}
            \mathbb{E}[\wbf_i(t)] = \zeros
            \quad \text{and} \quad
            \big\|\mathbb{E}[\wbf_i(t)\wbf_i(t)^T]\big\| \le C.
        \end{align*}

        \item\label{assum:st-var-growth}
        For all $\thetab \neq 0$ and all $t$,
        \begin{align}\label{eq:st-var-growth-z}
        \sup_{M, N} \frac{\thetab^T \cov\Big(\sum_{i = M+1}^{M+N} \zbf_i(t) \Big) \thetab}{\thetab^T \cov\Big(\sum_{i = 1}^N \zbf_i(t) \Big) \thetab} < \infty.
        \end{align}

        \item\label{assum:w-x-indep} For each $t$, all noises $\{\wbf_i(t)\}_{i = 1}^\infty$ are independent of contemporaneous states $\{\xbf_i(t)\}_{i = 1}^\infty$ at the same time $t$.

        \item\label{assum:f-bdd} The functions $f_i(.)$ are uniformly bounded, i.e., there exists $M < \infty$ such that for all $i, \xbf_{\Nc_i}, \wbf_i$
        \begin{align*}
            ||f_i(\xbf_{\Nc_i}, \wbf_i)|| \le M.
        \end{align*}

        \item\label{assum:rhomixing} For all $t$, the sequence of random vectors $\{\ybf_i(t)\}_{i = 1}^\infty$ in~\eqref{eq:rho-mixing-1d} is $\rho$-mixing with a slowly-varying residual factor $h(N)$.
    \end{enumerate}
\end{assumption}

The key assumption in Assumption~\ref{assum} is the $\rho$-mixing assumption in~\ref{assum:rhomixing} which, in essence, prevents global synchrony. In other words, it prevents a case where all $\{\ybf_i(t)\}_{i = 1}^\infty$ are strongly correlated with each other, in which case averaging does not do much---the average of $\{\ybf_i(t)\}_{i = 1}^\infty$ becomes very similar to any one of them.

The remaining 4 assumptions are technical and automatically satisfied in most practical applications. Assumptions~\ref{assum:w} and~\ref{assum:st-var-growth} are parallel to~\eqref{prop:cond_1} and~\eqref{prop:cond_3} in Proposition~\ref{prop:clt} and described thereafter. Next, assumption~\ref{assum:w-x-indep} asks for lack of instantaneous effect, where each noise variable takes at least one time step before affecting the state. Finally, note that Assumption~\ref{assum:f-bdd} is practically equivalent to uniform boundedness of \textit{state-noise trajectories}. If $(\xbf_{\Nc_i}, \wbf_i)$ are guaranteed to remain within a compact region $\Dc$, as is the case for any real-world state/noise variable, then $f_i$ needs to be bounded only on $\Dc$ and its value can be replaced by an arbitrary finite value (say, 0) outside of $\Dc$ without affecting the solutions of~\eqref{eq:ode}. This is true, e.g., if $f_i$ is any analytic function on $\Dc$.

With the stated assumptions above, we are ready to characterize the asymptotic spatial aggregate dynamics of the population of dynamical systems in~\eqref{eq:ode}, as follows.

\begin{theorem}\longthmtitle{Linearizing effect of spatial averaging on sequences of dynamical systems}\label{thm:s_avg}
Consider the population dynamics~\eqref{eq:ode}, and assume that Assumptions~\ref{assum:w}-\ref{assum:rhomixing} hold. Define the population average variables %
\begin{align}
    \mathbf{\bar x}(t) & = \frac{1}{\sqrt {Nh(N)}} \sum_{i = 1}^{N} \xbf_i(t) - \mathbb{E}[\xbf_i(t)], \label{eq:s_avg_state_full} \\
    \notag \mathbf{\bar w}(t) & = \frac{1}{\sqrt {Nh(N)}} \sum_{i = 1}^{N} \wbf_i(t),
\end{align}
and $\bar \ybf(t) = \begin{bmatrix} \bar \xbf(t)^T & \bar \wbf(t)^T \end{bmatrix}^T$ and assume that $\lim_{N \to \infty} \cov\big(\bar \xbf(t+1), \bar \ybf(t)\big)$ exists.
Then the relationship between $\mathbf{\bar x}(t+1)$, $\mathbf{\bar x}(t)$ and $\mathbf{\bar w}(t)$ becomes asymptotically linear as $N \to \infty$, i.e., for all $\xib \in \real^n, \omegab \in \real^m$,
\begin{align}\label{eq:cond-exp-s_avg}
    \notag \mathbb{E}\bigg [ \mathbf{\bar x}(t+1) \bigg | \begin{bmatrix}
    \mathbf{\bar x}(t) \\ \mathbf{\bar w}(t)
    \end{bmatrix} = \begin{bmatrix}
    \xib \\ \omegab
    \end{bmatrix} \bigg ]
    &\to \Abf_\infty(t)\xib  + \Bbf_\infty(t) \omegab
    \\
    \text{as} \ N &\to \infty,
\end{align}
where
\begin{align}\label{eq:s_AB}
\notag \Abf_\infty(t) &= \cov\big(\bar \xbf_\infty(t+1), \bar \xbf_\infty(t)\big) \cov\big(\bar \xbf_\infty(t)\big)^\dagger
\\
\notag \Bbf_\infty(t) &= \cov\big(\bar \xbf_\infty(t+1), \bar \wbf_\infty(t)\big) \cov\big(\bar \wbf_\infty(t)\big)^\dagger
\\
\notag \bar \xbf_\infty(t) &= \dlim_{N \to \infty} \bar \xbf(t)
\\
\bar \wbf_\infty(t) &= \dlim_{N \to \infty} \bar \wbf(t).
\end{align}
\end{theorem}
\begin{proof}
Fix $t$, and let $\zbf_i(t)$ be as in~\eqref{eq: conc_states}. From~\eqref{eq:ode}, Assumption~\ref{assum:rhomixing}, and Lemma~\ref{lem:rho_mixing_seq}, the sequence $\{\zbf_i(t)\}_{i = 1}^\infty$ is also $\rho$-mixing. Since mean-centering preserves correlation coefficients, $\{\zbf_i(t) - \E[\zbf_i(t)]\}_{i = 1}^\infty$ is $\rho$-mixing as well. The latter sequence is mean-zero by construction and finite-variance by Assumptions~\ref{assum:f-bdd} and~\ref{assum:w}. Furthermore, $h(N)$ is a residual factor for both $\{\xbf_i(t+1)\}_{i = 1}^\infty$ and $\big\{\begin{bmatrix} \xbf_i(t)^T & \wbf_i(t)^T \end{bmatrix}^T\big\}_{i = 1}^\infty$, so both of the limits
\begin{align*}
\lim_{N \to \infty} \cov\big(\bar \xbf(t+1)\big)
\quad \text{and} \quad
\lim_{N \to \infty} \cov\big(\bar \ybf(t)\big)
\end{align*}
exist. Since the limit of the covariance between $\bar \xbf(t+1)$ and $\bar \ybf(t)$ also exists by assumption, $h(N)$ is a (slowly-varying) residual factor for $\{\zbf_i(t) - \E[\zbf_i(t)]\}_{i = 1}^\infty$. These, together with Assumption~\ref{assum:st-var-growth}, satisfy all the requirements of Theorem~\ref{thm:multivar_clt_rho_mixing}. Therefore,
\begin{align*}%
    \mathbf{\bar z}(t) = \frac{1}{\sqrt{Nh(N)}} \sum_{i=1}^N \zbf_i(t) - \mathbb{E}[\zbf_i(t)] = \begin{bmatrix}
    \bar \xbf(t+1) \\ \bar \xbf(t) \\ \bar \wbf(t)
    \end{bmatrix}
\end{align*}
converges in distribution to $\bar \zbf_\infty(t) = \begin{bmatrix} \bar \xbf_\infty(t+1)^T & \bar \xbf_\infty(t)^T & \bar \wbf_\infty(t)^T \end{bmatrix}^T \sim \Nc(\zeros, \cov(\bar \zbf_\infty(t)))$, where $\cov(\bar \zbf_\infty(t)) = \lim_{N \to \infty} \cov(\bar \zbf(t))$. This, together with the fact that $\bar \zbf(t)$ is uniformly integrable (since it has finite variance), implies convergence of expectations over these distributions. In particular, for any $\xib, \omegab$,
\begin{align}\label{eq:cond-expec}
\notag \lim_{N \to \infty} &\E\bigg [ \mathbf{\bar x}(t+1) \bigg | \begin{bmatrix}
    \mathbf{\bar x}(t) \\ \mathbf{\bar w}(t)
    \end{bmatrix} = \begin{bmatrix}
    \xib \\ \omegab
    \end{bmatrix} \bigg ]
    \\
    &=
    \E\bigg [ \mathbf{\bar x_\infty}(t+1) \bigg | \begin{bmatrix}
    \mathbf{\bar x_\infty}(t) \\ \mathbf{\bar w_\infty}(t)
    \end{bmatrix} = \begin{bmatrix}
    \xib \\ \omegab
    \end{bmatrix} \bigg ].
\end{align}
Since $\zbf_\infty$ is normally distributed, the right hand side in~\eqref{eq:cond-expec} equals
\begin{align*}
\cov\bigg(\bar \xbf_\infty(t+1), \begin{bmatrix} \bar \xbf_\infty(t) \\ \bar \wbf_\infty(t) \end{bmatrix}\bigg) \cov\bigg(\begin{bmatrix} \bar \xbf_\infty(t) \\ \bar \wbf_\infty(t) \end{bmatrix}\bigg)^\dagger \begin{bmatrix} \xib \\ \omegab \end{bmatrix}
\end{align*}
which simplifies to $\Abf_\infty(t) \xib + \Bbf_\infty(t) \omegab$ in~\eqref{eq:cond-exp-s_avg} because of the independence of $\bar \xbf_\infty(t)$ and $\bar \wbf_\infty(t)$ resulting from Assumption~\ref{assum:w-x-indep}.
\end{proof}

As far as linearity of the limit dynamics is concerned, it is important to note that the matrices $\Abf_\infty(t)$ and $\Bbf_\infty(t)$ in Theorem~\ref{thm:s_avg} depend on the \emph{distributions} of states $\{\xbf_i(t)\}_{i = 1}^\infty, \{\xbf_i(t+1)\}_{i = 1}^\infty$ and noises $\{\wbf_i(t)\}_{i = 1}^\infty$, but not on their specific values (realizations). In other words, for any choice of model in~\eqref{eq:ode}, i.e., for any choice of $\{f_i(\cdot)\}_{i = 1}^\infty$, $p_0$, and $p_w$, the matrices $\Abf_\infty(t)$ and $\Bbf_\infty(t)$ can be in theory pre-computed for all $t \ge 0$, and would be the same for all realizations of noises and states generated by this system. This is important for the dynamics to be actually linear, as illustrated by the following simple example.

\begin{example}\longthmtitle{LTV vs. nonlinear dynamics}\label{ex:ei}
  Consider a dynamical system described by
  \begin{align*}
      y(t+1) = a(t) y(t), \quad y(0) \sim \Nc(e, 1)
  \end{align*}
  and three cases, as follows. \\
  \begin{enumerate}[wide]
      \item $a(t) = e^{2^t} $. Here the system is clearly LTV. \\
      \item $a(t) = y(t) $. Here the system is clearly nonlinear. \\
      \item $a(t) = \mathbb{E}[y(t)] $. Although the system here may look nonlinear, the coefficient $\mathbb{E}[y(t)]$ can be pre-computed for all $t \ge 0$ based only on the distribution of $y(0)$. In fact, it is straightforward to see that here $a(t) = e^{2^t}$, making this system equivalent to that in~(i). \oprocend
  \end{enumerate}

\end{example}

\subsection{Time-Invariant Limit Dynamics}\label{subsec:lti}

Interestingly, if the dynamics in~\eqref{eq:ode} satisfy additional assumptions, the average dynamics tend not only to a linear system but further to an LTI one. For ease of notation, let
\begin{align}\label{eq:pop-dynamics}
    \xbf(t+1) = F(\xbf(t), \wbf(t))
\end{align}
denote the combined dynamics of all sub-systems in~\eqref{eq:ode},
where $\xbf(t) \in \real^{Nn}$ and $\wbf(t) \in \real^{Nm}$ are the concatenations of all states and noises, respectively. We also need to make additional assumptions, as follows.

\begin{assumption}\longthmtitle{Additional assumptions for time-invariance of limit dynamics}\label{assum2}
\begin{enumerate}[label=\rm(A\arabic*), leftmargin=25pt]
\setcounter{enumi}{5}
\item\label{assum:pwstat} The noise distribution $p_w(t) \equiv p_w$ is time-invariant.
\item\label{assum:global-Lipschitz} The function $F(\cdot)$ is globally Lipschitz in $\xbf$, i.e., for any $\wbf \in \real^{Nm}$ there exists $L(\wbf) \ge 0$ such that
    \begin{align}\label{lips}
    \|F(\xbf_1, \wbf) - F(\xbf_2, \wbf)\| \le L(\wbf) \|\xbf_1 - \xbf_2\|,
\end{align}
for all $\xbf_1, \xbf_2 \in \real^{Nn}$.
\item\label{assum:p0} the initial distribution $p_0$ is such that $p_0(A) = 0$ for any set $A$ where $p^{*}(A) = 0$;
\item\label{assum:pw} for any set $A$ with $p^{*}(A) > 0$, the noise distribution $p_w$ is such that for all $\xbf \in \real^{Nn}$, $\P\{F(\xbf, \wbf) \in A\} > 0$. %
\item\label{assum:uniformly-rho-mixing} The sequence $\{\ybf_i(t)\}_{i = 1}^\infty$
is $\rho$-mixing \emph{uniformly across time}, i.e., $\sup_t \varrho(n, t) \to 0$ as $n \to \infty$, where $\varrho(n, t)$ is defined as in~\eqref{eq:rhoto0} for the sequence $\{\ybf_i(t)\}_{i = 1}^\infty$. %
\oprocend
\end{enumerate}
\end{assumption}

Assumption~\ref{assum:pwstat} is clearly necessary if we expect~\eqref{eq:ode} to admit a stationary solution. %
Assumption~\ref{assum:global-Lipschitz}, while being restrictive on the space of all functions, becomes mild when considering only bounded functions as required by assumption~\ref{assum:f-bdd}. %
Assumptions~\ref{assum:p0} and~\ref{assum:pw} are technical and ensure the distributions $p_0$ and $p_w$ are well-behaved. Assumption~\ref{assum:uniformly-rho-mixing} is also technical and prevents contrived cases where, e.g., the system gradually converges to a state of global synchrony.

The next result uses Assumptions~\ref{assum:pwstat}-\ref{assum:pw} to show the existence and attractivity of stationary solutions, which will then be used to prove the LTI version of Theorem~\ref{thm:s_avg}.

\begin{theorem}\longthmtitle{Existence and attractivity of stationary solutions}\label{thm:stationary_sol}
Consider the population dynamics in~\eqref{eq:ode} and assume that Assumptions~\ref{assum:w}-\ref{assum:global-Lipschitz} hold.
Then the population has a stationary solution
\begin{align}\label{eq:stat-sol}
    \xbf^{*}(t) \sim p^{*}, \qquad t \ge 0.
\end{align}
where $p^*$ is independent of time. If, further, Assumptions~\ref{assum:p0}-\ref{assum:pw} hold, then $\xbf(t)$ converge to $\xbf^{*}(t)$ in distribution as $t \to \infty$.
\end{theorem}
\begin{proof}
    The proof is the same as proof of~\cite[Thm IV.4]{ahmed2022linearizing}.
\end{proof}

Combining the linearity of Theorem~\ref{thm:s_avg} and the stationarity of Theorem~\ref{thm:stationary_sol} ensures the convergence of the average population dynamics to an LTI system, as shown next.

\begin{theorem}\longthmtitle{LTI average population dynamics}
\label{thm:lti}
Consider the population dynamics~\eqref{eq:ode} and assume that assumptions~\ref{assum:w}-\ref{assum:uniformly-rho-mixing} hold.
Then,
\begin{align*}
    \mathbb{E}\bigg [ \mathbf{\bar x}(t+1) \bigg | \begin{bmatrix}
    \mathbf{\bar x}(t) \\ \mathbf{\bar w}(t)
    \end{bmatrix} = \begin{bmatrix}
    \upsilonb_0 \\ \zetab_0
    \end{bmatrix} \bigg ]
    \to \Abf_\infty^* \upsilonb_0 + \Bbf_\infty^* \zetab_0
\end{align*}
as $N, t \to \infty$, where
\begin{align}\label{eq:limab}
\Abf_\infty^* = \lim_{t \to \infty} \Abf_\infty(t), \quad \Bbf_\infty^* = \lim_{t \to \infty} \Bbf_\infty(t),
\end{align}
and $\Abf_\infty(t), \Bbf_\infty(t)$ are as in~\eqref{eq:s_AB}.
\end{theorem}
\begin{proof}
By Theorem~\ref{thm:stationary_sol} %
each subsystem $\xbf_i(t)$ converges in distribution to the corresponding stationary solution $\xbf_i^{*} \sim p_i^{*}$ as $t \to \infty$, which are $\rho$-mixing (in the limit) by~\ref{assum:uniformly-rho-mixing}. %
Define
\begin{align*}
    \bar \xbf^{*} &= \frac{1}{\sqrt {Nh(N)}} \sum_{i = 1}^{N} \xbf_i^{*} - \mathbb{E}[\xbf_i^{*}],
\end{align*}
and let $\bar \xbf^{*}_\infty \sim \Nc(\zeros, \Sigmab_{\bar x^*})$ where $\Sigmab_{\bar x^*} = \cov(\bar \xbf^*)$. Then by %
Theorem~\ref{thm:multivar_clt_rho_mixing}, $\bar \xbf^*$ converges in distribution to $\bar \xbf^*_\infty$ as $N \to \infty$ and, therefore,
\begin{align}\label{eq:x_star}
    \dlim_{N,t \to \infty} %
    \bar \xbf(t) = \dlim_{t \to \infty} \bar \xbf_\infty(t) = \bar \xbf^{*}_\infty.
\end{align}
It then follows from~\eqref{eq:x_star} that all the time-dependent covariances in~\eqref{eq:s_AB} converge to their respective limits and, hence, so do the matrices $\Abf_\infty(t)$ and $\Bbf_\infty(t)$. This completes the proof.
\end{proof}

Similar to the guaranteed linearity of Theorem~\ref{thm:s_avg}, the result of Theorem~\ref{thm:lti} is asymptotic and LTI dynamics are only \textit{approached} as $N$ and $t$ grow to infinity. In many real-world systems, however, dynamics settle to stationary solutions after only a few time steps. Further, and more importantly for our discussion, nonlinearities often vanish rather quickly with averaging. We next focus on formalizing the latter dimension, namely, the rate of convergence to linearity under spatial averaging.

\subsection{Rate of Convergence to Linearity}\label{sub_sec:finite_err}

In this section we seek to estimate the rate at which spatially-averaged dynamics converge to linearity. We have empirically shown that this rate is often rather fast~\cite{ahmed2022linearizing}, but rigorous characterizations are missing. The main result of this section is the following, where we prove that convergence of conditional expectations (similar to~\eqref{eq:cond-exp-s_avg}) occurs at the $O(1/\sqrt N)$ rate. We present this result for the simplest  case of i.i.d. sequences of static random vectors and discuss generalizations afterwards.

\begin{theorem}\longthmtitle{Rate of convergence to linearity}\label{thm:rate}
Consider a sequence of $N$ i.i.d random vectors $\zbf_1, \zbf_2, \dots, \zbf_N \in \real^{n+q}$ where, for all $i$,
\begin{align}\label{eq:rate-model}
\zbf_i &= \begin{bmatrix} \xbf_i \\ \ybf_i \end{bmatrix}, \quad \xbf_i \in \real^n, \quad \ybf_i \in \real^q
\\
\notag \E[\zbf_i] &= \mub_z, \quad \cov(\zbf_i) = \Sigmab_z, \quad \|\zbf_i\|_\infty \le M < \infty.
\end{align}
Let
\begin{align*}
    \bar \zbf = \frac{1}{N} \sum_{i=1}^N \zbf_i \qquad \text{and} \qquad \zbf^* = \begin{bmatrix} \xbf^* \\ \ybf^* \end{bmatrix} \sim \Nc\Big(\mub_z, \frac{1}{N} \Sigmab_z\Big)
\end{align*}
Then, for any $\upsilonb \in \real^q$ and any $\delta > 0$,
\begin{align}\label{eq: mean_cond_dis_diff}
    \Big\|\E[\bar \xbf| \bar \ybf \in \Bc_{\delta}(\upsilonb)] - \E[\xbf^* | \ybf^* \in \Bc_{\delta}(\upsilonb)] \Big\| = O\Big(\frac{1}{\sqrt N}\Big),
\end{align}
where $\Bc_\delta(\upsilonb) = \{\ybf \in \real^q:  ||\ybf - \upsilonb|| \le \delta \}$ is a $\delta$-ball around $\upsilonb$.
\end{theorem}
\begin{proof}
According to Berry-Essen theorem~\cite{raivc2019multivariate}, for any convex sets $U \subseteq \real^{n+q}$,
\begin{align}\label{eq:berry-esseen}
    \big|\mathbb{P}\big(\bar \zbf \in U \big) - \mathbb{P}\big(\zbf^* \in U \big)\big| \le \frac{C}{\sqrt N}
\end{align}
where $C$ is a constant independent of $N$. So consider, in particular, the sets of the form
\begin{align*}
U = \big\{\xbf \in \real^n | x_i \le \xi \big\} \times \Bc_\delta(\upsilonb)
\end{align*}
for any $i \in \{1, \dots, n\}$ and any $\xi \in \real$. Plugging this into~\eqref{eq:berry-esseen} we get
\begin{align}\label{eq:cdf_joint}
    \big|\P\big\{\bar x_i \le \xi, \bar \ybf \in \Bc_\delta(\upsilonb)\big\} -\P\big\{x^*_i \le \xi, \ybf^* \in \Bc_\delta(\upsilonb)\big\}\big| \le \frac{C}{\sqrt N}.
\end{align}
Similarly, by plugging $U = \real^n \times \Bc_\delta(\upsilonb)$ in~\eqref{eq:berry-esseen} we get
\begin{align}\label{eq:cdf_cond_var}
    \big|\P\big\{\bar \ybf \in \Bc_\delta(\upsilonb)\big\} - \P\big\{\ybf^* \in \Bc_\delta(\upsilonb)\big\}\big| \le \frac{C}{\sqrt N}.
\end{align}
Further, let
\begin{align}\label{eq:cdf_lower_bound}
    \epsilon = \min\Big\{\P\big\{\bar \ybf \in \Bc_\delta(\upsilonb)\big\} \ , \ \P\big\{\ybf^* \in \Bc_\delta(\upsilonb)\big\}\Big\} > 0.
\end{align}
Combining~\eqref{eq:cdf_joint}, \eqref{eq:cdf_cond_var} and~\eqref{eq:cdf_lower_bound} we get
\begin{align}\label{eq:diff-cdfs}
    \notag &\Big|F_{\bar x_i | \bar \ybf}\big(\xi | \Bc_\delta(\upsilonb)\big) - F_{x^*_i | \ybf^*}\big(\xi | \Bc_\delta(\upsilonb)\big)\Big|
    \\
    \notag  &\triangleq \Big|\P\big\{\bar x_i \le \xi \ | \ \bar \ybf \in \Bc_\delta(\upsilonb)\big\} - \P\big\{x^*_i \le \xi \ | \ \ybf^* \in \Bc_\delta(\upsilonb)\big\}\Big|
    \\
    \notag  &= \bigg | \frac{\P\big\{\bar x_i \le \xi \ , \ \bar \ybf \in \Bc_\delta(\upsilonb)\big\}}{\P\big\{\bar \ybf \in \Bc_\delta(\upsilonb)\big\}} - \frac{\P\big\{x^*_i \le \xi \ , \ \ybf^* \in \Bc_\delta(\upsilonb)\big\}}{\P\big\{\ybf^* \in \Bc_\delta(\upsilonb)\big\}} \bigg |
    \\
     &\le \frac{2C}{\sqrt N \epsilon^2}
\end{align}
where in the latter inequality we used $\big|\frac{a}{b} - \frac{c}{d}\big| \le \frac{|d| |a-c| + |c| |b-d|}{|b||d|}$.
Now, using the standard relationship between $\E[X] = -\int_0^\infty F_X(x) dx + \int_0^\infty (1-F_X(x))dx$ between the cumulative distribution function (CDF) and expected value of any random variable, we can translate the bound on the difference of conditional CDFs in~\eqref{eq:diff-cdfs} into a bound on the difference of conditional expectations, as follows.
\begin{align}\label{eq:sum-of-ints}
    \notag & \Big| \E\big[\bar x_i \ | \ \bar \ybf \in \Bc_\delta(\upsilonb)\big] -  \E\big[x^*_i \ | \  \ybf^* \in \Bc_\delta(\upsilonb)\big] \Big|
    \\
    \notag &= \bigg| \int_{-\infty}^{0} \Big(F_{x^*_i | \ybf^*}\big(\xi | \Bc_\delta(\upsilonb)\big) - F_{\bar x_i | \bar \ybf}\big(\xi | \Bc_\delta(\upsilonb)\big)\Big) d\xi
    \\
    \notag &\qquad+ \int_0^\infty \Big(F_{x^*_i | \ybf^*}\big(\xi | \Bc_\delta(\upsilonb)\big) - F_{\bar x_i | \bar \ybf}\big(\xi | \Bc_\delta(\upsilonb)\big)\Big) d\xi  \bigg|
    \\
    \notag &\le  \int_{-\infty}^{0} \Big|F_{x^*_i | \ybf^*}\big(\xi | \Bc_\delta(\upsilonb)\big) - F_{\bar x_i | \bar \ybf}\big(\xi | \Bc_\delta(\upsilonb)\big)\Big| d\xi
    \\
    &\qquad+ \int_0^\infty \Big|F_{x^*_i | \ybf^*}\big(\xi | \Bc_\delta(\upsilonb)\big) - F_{\bar x_i | \bar \ybf}\big(\xi | \Bc_\delta(\upsilonb)\big)\Big| d\xi
\end{align}
Note that, because the bound in~\eqref{eq:diff-cdfs} is independent of $\xi$, substituting it directly into~\eqref{eq:sum-of-ints} will give a useless bound of $\infty$. However, because $x_i \in [-M, M]$ by assumption, $F_{\bar x_i | \bar \ybf}\big(\xi | \Bc_\delta(\upsilonb)\big) = 0$ for all $\xi \notin [-M, M]$ and we can use this to achieve a finite bound, as follows. For reasons that we will see later, let $M' \ge M$ be a (yet-to-be-determined) more conservative bound on $x_i$. Then,
\begin{align}\label{eq:bound-plus-int}
    \notag &\Big| \E\big[\bar x_i \ | \ \bar \ybf \in \Bc_\delta(\upsilonb)\big] -  \E\big[x^*_i \ | \  \ybf^* \in \Bc_\delta(\upsilonb)\big] \Big|
    \\
    \notag &\le \int_{-\infty}^{-M'} F_{x^*_i | \ybf^*}\big(\xi | \Bc_\delta(\upsilonb)\big) d\xi + \int_{-M'}^{0}  \frac{2C}{\sqrt N \epsilon^2} d\xi
    \\
    \notag &\qquad+ \int_{0}^{M'} \frac{2C}{\sqrt N \epsilon^2} d\xi +
    \int_{M'}^{\infty} \Big ( 1 - F_{x^*_i | \ybf^*}\big(\xi | \Bc_\delta(\upsilonb)\big)  \Big) d\xi \\
    &= \frac{4M'C}{\sqrt N \epsilon^2} + 2 \int_{-\infty}^{-M'} F_{x^*_i | \ybf^*}\big(\xi | \Bc_\delta(\upsilonb)\big) d\xi.
\end{align}
To evaluate the integral in~\eqref{eq:bound-plus-int}, note that
\begin{align*}
&F_{x^*_i | \ybf^*}\big(\xi | \Bc_\delta(\upsilonb)\big) = \P\big\{x^*_i \le \xi \ | \ \ybf^* \in \Bc_\delta(\upsilonb)\big\}
\\
&\le \sup_{\ybf_0 \in \Bc_\delta(\upsilonb)} \P\big\{x^*_i \le \xi \ | \ \ybf^*  = \ybf_0 \big\}
\\
&= \sup_{\ybf_0 \in \Bc_\delta(\upsilonb)} \Phi\Bigg( \frac{\xi - \mu_{x_i} - \Sigmab_{x_i \ybf} \Sigmab_\ybf^{-1} (\ybf_0 - \mub_\ybf)}{\big[\frac{1}{N}\big(\Sigma_{x_i} - \Sigmab_{x_i \ybf} \Sigmab_\ybf^{-1} \Sigmab_{\ybf x_i}\big)\big]^{-\frac{1}{2}}} \Bigg)
\\
&= \Phi\Bigg( \frac{\xi - \mu_{x_i} - \inf_{\ybf_0 \in \Bc_\delta(\upsilonb)} \Sigmab_{x_i \ybf} \Sigmab_\ybf^{-1} (\ybf_0 - \mub_\ybf)}{\big[\frac{1}{N}\big(\Sigma_{x_i} - \Sigmab_{x_i \ybf} \Sigmab_\ybf^{-1} \Sigmab_{\ybf x_i}\big)\big]^{\frac{1}{2}}} \Bigg)
\\
&= \Phi\Bigg( \frac{\xi - \mu_{x_i} - \Sigmab_{x_i \ybf} \Sigmab_\ybf^{-1} (\upsilonb - \mub_\ybf) + \delta \|\Sigmab_{x_i \ybf} \Sigmab_\ybf^{-1}\|}{\big[\frac{1}{N}\big(\Sigma_{x_i} - \Sigmab_{x_i \ybf} \Sigmab_\ybf^{-1} \Sigmab_{\ybf x_i}\big)\big]^{\frac{1}{2}}} \Bigg)
\end{align*}
where $\Phi(\cdot)$ is the standard normal CDF. For ease of notation let $\xi_0 = \mu_{x_i} + \Sigmab_{x_i \ybf} \Sigmab_\ybf^{-1} (\upsilonb - \mub_\ybf) - \delta \|\Sigmab_{x_i \ybf} \Sigmab_\ybf^{-1}\|$ and $\sigma_0^2 = \Sigma_{x_i} - \Sigmab_{x_i \ybf} \Sigmab_\ybf^{-1} \Sigmab_{\ybf x_i}$. Then
\begin{align*}
F_{x^*_i | \ybf^*}\big(\xi | \Bc_\delta(\upsilonb)\big) \le \Phi\Big(\frac{\xi \!-\! \xi_0}{\sigma_0} \sqrt N \Big)
\!=\! \frac{\sqrt N}{\sqrt{2\pi} \sigma_0} \int_0^\xi \!\!e^{-\frac{(\zeta - \xi_0)^2}{2 \sigma_0^2} N} \!\!d\zeta
\end{align*}
Notice that, from~\eqref{eq:bound-plus-int}, always $\xi \le -M'$, where $M'\ge M$ is yet to be determined. Choose $M' = \max\{M, 1 - \xi_0\}$. Then, for all $\xi \le -M'$ we have
\begin{align*}
F_{x^*_i | \ybf^*}\big(\xi | \Bc_\delta(\upsilonb)\big) \le \frac{\sqrt N}{\sqrt{2\pi} \sigma_0} \int_0^\xi e^{-\frac{|\zeta - \xi_0|}{2 \sigma_0^2} N} d\zeta
= \frac{\sqrt 2 \sigma_0}{\sqrt{\pi N}} e^{\frac{\xi - \xi_0}{2 \sigma_0^2} N}
\end{align*}
Substituting this into~\eqref{eq:bound-plus-int} and integrating from $-\infty$ to $-M'$ gives
\begin{align*}
&\Big| \E\big[\bar x_i \ | \ \bar \ybf \in \Bc_\delta(\upsilonb)\big] -  \E\big[x^*_i \ | \  \ybf^* \in \Bc_\delta(\upsilonb)\big] \Big|
\\
&\le \frac{4M'C}{\sqrt N \epsilon^2} + \frac{2}{\sqrt \pi} \Big(\frac{2 \sigma_0}{N}\Big)^{\frac{3}{2}} e^{-\frac{M'+\xi_0}{2\sigma_0^2}N} = O\Big(\frac{1}{\sqrt N}\Big).
\end{align*}
Because this is true for all $i = 1, \dots, n$, and $\|\cdot\| \le \sqrt n \|\cdot\|_\infty$, we get \eqref{eq: mean_cond_dis_diff}, completing the proof.
\end{proof}

A number of remarks about Theorem~\ref{thm:rate} are in order. First, this is the only result in this paper where we assume independence among the subsystems, in order to simplify the analysis. If the subsystems form a $\rho$-mixing sequence, similar results can be obtained using extensions of the Berry-Esseen theorem (see, e.g.,~\cite{wang2019berry}) but the analysis becomes more involved and the rates of convergence become variable depending on the rate at which correlations decay (e.g., $n^{-\frac{1}{6}} \log n$ if correlations decay as $n^{-\frac{3}{2}}$). Second, the $L_\infty$ nature of the Berry-Esseen theorem (i.e., that the bound in~\eqref{eq:berry-esseen} does not depend on the set $U$) has forced us to move from the singleton-conditioning in~\eqref{eq:cond-exp-s_avg} to neighborhood-conditioning in~\eqref{eq: mean_cond_dis_diff}. To the best of our knowledge this is not avoidable for continuous-valued conditioning variables (see~\cite{dey2023stein} for conditioning on discrete variables) but the radius of $\Bc_\delta(\upsilonb)$ can be made sufficiently small so that conditioning on $\bar \ybf \in \Bc_\delta(\upsilonb)$ is a good approximation of conditioning on $\bar \ybf = \upsilonb$. Finally, we presented Theorem~\ref{thm:rate} for static random variables, but its extension to dynamical systems is straightforward, with $\xbf_i$ and $\ybf_i$ in~\eqref{eq:rate-model} being replaced by $\xbf_i(t+1)$ and $\ybf_i(t) = \begin{bmatrix} \xbf_i(t)^T & \wbf_i(t)^T \end{bmatrix}^T$, respectively. Under the same assumptions, namely, i.i.d. subsystems with uniformly bounded states and noises, the same conclusion as in~\eqref{eq: mean_cond_dis_diff} would then hold for conditional expectation of average states at any time $t+1$ given neighborhood (not singleton) conditions on average states and noises at the previous time $t$.%

\section{Linearizing Effect of Spatial Averaging on Spatially-Embedded Dynamical Systems}\label{sec:3d}

\subsection{Problem Formulation}

While the results of Section~\ref{sec:s_avg_1d} are general in many respects (form of nonlinearity, noise characteristics, etc.) they are still limited in requiring a decay of correlations with the difference in two subsystems' linear indices (i.e., $\varrho(|i - j|) \to 0$ as $|i - j| \to \infty$). In brain networks that initially motivated the present work~\cite{nozari2023macroscopic}, e.g., it is well-known that correlations between neurons decay as their distance grows~\cite{smith2008spatial,smith2013spatial,rosenbaum2017spatial}, but this distance is physical (Euclidean) rather than between linear indices. Similarly, common forms of spatial averaging that occur in the brain, such as those underlying functional magnetic resonance imaging (fMRI) or electroencephalography (EEG), occur over neurons that are close to each other in a physical sense rather than index-wise. Therefore, in this section we generalize our results of Section~\ref{sec:s_avg_1d} to systems where both the decay of pairwise correlations and the domains of spatial averaging occur in Euclidean space.

Consider the same form of nonlinear population dynamics~\eqref{eq:ode} but now assume that each subsystem $i$ lies at a location $\rbf_i \in \real^d$. These locations could be deterministic and fixed, or themselves be stochastic, and we will return to this later. For the sake of generality we keep $d$ arbitrary, though we are often interested in $d = 3$. Instead of~\eqref{eq:x_Ni} let
\begin{align}\label{eq:x_Ni_3d}
    \xbf_{\Nc_i}(t) =
    \begin{bmatrix}
        \xbf_j(t) \ : \ \|\rbf_i - \rbf_j\| \le \tau
    \end{bmatrix},
\end{align}
where $\tau > 0$ can still be arbitrarily large but finite. We can now formulate the problem that we will tackle in this section, which closely parallels Problem~\ref{prob1}.

\begin{problem}\longthmtitle{Linearizing Effect of Spatial Averaging on Spatially-Embedded Dynamical Systems}\label{prob2}
    Consider a heterogeneous population of nonlinear dynamical systems described by~\eqref{eq:ode} and~\eqref{eq:x_Ni_3d}. Define the population's average state vector as
    \begin{align}\label{eq:s_avg_state_3d}
        \bar \xbf(t) = \frac{1}{\phi(N_\Rc)} \sum_{i: \|\rbf_i\| \le \Rc} \xbf_i(t) - \mathbb{E}[\xbf_i(t)],
    \end{align}
    where, for any $\Rc > 0$,
    \begin{align}\label{eq:NR}
    N_\Rc = \big|\{i \ | \ \|\rbf_i\| \le \Rc\}\big|,
    \end{align}
    and $\phi(\cdot)$ is a normalization factor as in~\eqref{eq:s_avg_state}. Prove, under appropriate assumptions and choice of $\phi(\cdot)$, that
    the dynamics of $\bar \xbf(t)$ becomes asymptotically linear as $\Rc \to \infty$. \oprocend
\end{problem}

The rest of this section is concerned with addressing Problem~\ref{prob2}. We will follow the same general steps as in Section~\ref{sec:s_avg_1d}, focusing on what needs to be done differently in order to address the challenges that arise from the new formulation.

\subsection{$\rho^*$-Mixing Sequences}\label{sub_sec:rhostar-mixing-seq}

Consider a spatially-embedded sequence of dynamical systems as in Problem~\ref{prob2}, which gradually grows in size as $\Rc, N_\Rc \to \infty$. Regardless of how the subsystems are enumerated, if all we know is that their correlations decay with their physical distance, there is no guarantee that the same happens if $|i-j| \to \infty$. Therefore, to address Problem~\ref{prob2} we need a different notion of mixing, as introduced next.

\begin{definition}\longthmtitle{Multivariate $\rho^*$-mixing sequence}\label{def:rho_star_mixing}
Consider a sequence of random vectors $\xib_1, \xib_2, \dots \in \real^q$ with corresponding Euclidean positions $\rbf_1, \rbf_2, \dots \in \real^d$. For any $r > 0$,
define
\begin{align}\label{eq:rho_star_d}
    \!\!\varrho^*(r) = \sup_{A, B} \ \rho\Big(\sigma\big(\{\xib_i \ | \ \rbf_i \in A\}\big), \sigma\big(\{\xib_j \ | \ \rbf_j \in B\}\big)\Big),
\end{align}
where the supremum is taken over all nonempty sets $A, B \subset \real^d$ such that
\begin{align*}
    \inf_{\rbf_i \in A, \ \rbf_j \in B} \|\rbf_i - \rbf_j\| \ge r,
\end{align*}
and $\rho(\cdot, \cdot)$ is defined in~\eqref{eq:rhoAB}.
The sequence is $\rho^*$-mixing if
\begin{align}\label{eq: rho_d_to_zero}
    \varrho^*(r) \to 0 \quad \text{as} \quad r \to \infty.
\end{align}
\end{definition}

$\rho^*$-mixing sequences share many of the same properties with $\rho$-mixing ones. For example, the notion of a residual factor (Definition~\ref{def:resfact}) extends naturally to $\rho^*$-mixing sequences. The formal definition, for clarity, is as follows.

\begin{definition}\longthmtitle{Residual factor for spatially-embedded sequences}\label{def:resfact-star}
Let $\xib_1, \xib_2, \dots$ be a $\rho^*$-mixing sequence with corresponding Euclidean positions $\rbf_1, \rbf_2, \dots$. For any $\Rc > 0$ let
\begin{align*}
S_\Rc = \sum_{i: \|\rbf_i\| \le \Rc} \xib_i,
\end{align*}
and $N_\Rc$ be as in~\eqref{eq:NR}.
The function $h(N_\Rc)$ is called a \emph{residual factor} for $\{\xib_i\}_{i = 1}^\infty$ if
\begin{align*}
\lim_{\Rc \to \infty} \frac{1}{N_\Rc h(N_\Rc)} \cov(S_\Rc) < \infty,
\end{align*}
i.e., the limit exists and is finite. \oprocend
\end{definition}

In what follows we also need the following result, which parallels Lemma~\ref{lem:rho_mixing_seq} and can be proved using the same approach.

\begin{lemma}\longthmtitle{Sequence formed through transformation of $\rho^*$-mixing sequence is $\rho^*$-mixing}\label{lem:rho_star_mixing_seq}
Consider a sequence of random vectors $\xib_1, \xib_2, \dots \in \real^q$ with corresponding Euclidean positions $\rbf_1, \rbf_2, \dots \in \real^d$. For each $i$, let $\xib_{\Nc_i} \in \real^{\nu_i}$ be as in~\eqref{eq:x_Ni_3d}, $h_i$ be a measurable function from $\real^{\nu_i}$ to $\real^k$, and
\begin{align*}
    \zetab_i = h(\xib_{\Nc_i}) \in \real^k.
\end{align*}
If the sequence $\xib_1, \xib_2, \dots$ is $\rho^*$-mixing, then so is the sequence $\zetab_1, \zetab_2, \dots$. \oprocend
\end{lemma}

Next we move to the CLT, which lies at the core of our methodology. Here the main challenge lies in proving a result similar to Proposition~\ref{prop:clt} for $\rho^*$-mixing sequences, as tackled next. Once that is shown, generalization to multivariate sequences will be straightforward. %

\begin{theorem}\longthmtitle{CLT for scalar $\rho^*$-mixing sequences}\label{thm:clt_rho_star}
Consider a scalar $\rho^*$-mixing sequence of random variables $\xi_1, \xi_2, \dots \in \real$ with corresponding Euclidean positions $\rbf_1, \rbf_2, \dots \in \real^d$. %
For any $\Rc > 0$, let
\begin{align}\label{eq:slowly_varying_3d}
    S_\Rc = \sum_{i: \|\rbf_i\| \le \Rc} \xi_i, \qquad \sigma_\Rc^2 = \var(S_\Rc) = N_\Rc h(N_\Rc),
\end{align}
where $N_\Rc$ is defined in~\eqref{eq:NR} %
and $h(N_\Rc)$ is a slowly-varying function. Assume, further, that
\begin{subequations}
\begin{align}
    \label{eq:zero_mean_3d} &\!\!\!\!\mathbb{E}[\xi_i] = 0 \ \text{and} \ |\xi_i| \le M %
    \ \text{for all} \ i \ \text{and some} \ M < \infty
    \\
    \label{eq:unif_int_3d} &\!\!\!\!\!\!\sup_{\Rc > 0} \frac{\E\big[|S_{\gamma \Rc}|^p\big]}{\sigma_\Rc^p}  < \infty \ \text{for all} \ \gamma \in [0, 1] \ \text{and some} \ p > 2
    \\
    \label{eq:Ndinf} &\!\!\!\!N_\Rc \propto \Rc^d, \ \text{i.e.,} \ \lim_{\Rc \to \infty} \frac{N_\Rc}{\Rc^d} = C_v \ \text{for some} \ C_v > 0. %
\end{align}
\end{subequations}
Then
\begin{align}\label{eq:rho-star-mixing-clt}
    \frac{S_\Rc}{\sigma_\Rc} \stackrel{d}{\to} \Nc(0, 1) \quad \text{as} \quad \Rc \to \infty.
\end{align}
\end{theorem}
\begin{proof}
Consider the partial-sum stochastic process
\begin{align*}
    W_\Rc(\gamma) = \frac{S_{\gamma^\frac{1}{d} \Rc}}{\sigma_\Rc}, \quad \text{for all} \ \gamma \in [0, 1]. %
\end{align*}
In what follows we will prove the weak convergence of
$W_\Rc(\gamma)$ %
to the standard Brownian process $W(\gamma)$, %
from which the statement of the theorem follows as a special case.
By~\eqref{eq:zero_mean_3d} we have
\begin{align}\label{eq:mean-3d-0}
    \E[W_\Rc(\gamma)] = 0,
\end{align}
and by~\eqref{eq:slowly_varying_3d} and~\eqref{eq:Ndinf} we get
\begin{align}\label{eq:var-gamma}
    \notag &\lim_{\Rc \to \infty}  \E\big[W_\Rc^2(\gamma)\big] = \lim_{\Rc \to \infty}  \frac{1}{\sigma_{\Rc}^2} \E \Big[S_{\gamma^\frac{1}{d}\Rc}^2\Big]
    \\
    \notag &= \lim_{\Rc \to \infty}  \frac{N_{\gamma^\frac{1}{d} \Rc} h\Big(N_{\gamma^\frac{1}{d} \Rc}\Big)}{N_\Rc h(N_\Rc)}
    \\
    &= \lim_{\Rc \to \infty}  \frac{N_{\gamma^\frac{1}{d} \Rc}}{\gamma \Rc^d} \cdot \gamma \cdot \frac{\Rc^d}{N_\Rc} \cdot \frac{h\Big(N_{\gamma^\frac{1}{d} \Rc}\Big)}{h(N_\Rc)}
    = \gamma,
\end{align}
where in the last equality we used the fact that $\lim_{\Rc \to \infty} h\big(N_{\gamma^{1/d} \Rc}\big) / h(N_\Rc) = 1$ by the uniform convergence of slowly-varying functions~\cite[Thm. 1.2.1]{bingham1989regular}.

Further, by~\eqref{eq:unif_int_3d}
the collection of random variables,
\begin{align}\label{eq:uniint}
\big\{W_\Rc^2(\gamma) \ | \ \Rc \in \real \big\}
\end{align}
is uniformly integrable. %
To prove the weak convergence of $W_\Rc(\gamma)$ to $W(\gamma)$, we need to show that $W_\Rc(\gamma)$ satisfies two more conditions. The first is tightness, namely, that for any $\epsilon, \eta > 0$ there exists $ \delta > 0$ such that for all $|s - t| < \delta$ and all sufficiently large $\Rc$,
\begin{align}\label{eq: tightness}
    \mathbb{P}\Big( |W_\Rc(s) - W_\Rc(t)| \ge \epsilon \Big) \le \eta.
\end{align}
To show this, note that by Markov's inequality,
\begin{align}\label{eq:markov}
    \mathbb{P}\Big( |W_\Rc(s) - W_\Rc(t)| \ge \epsilon \Big) \le \frac{\mathbb{E}\big[|W_\Rc(s) - W_\Rc(t)|\big]}{\epsilon}.
\end{align}
Further, from~\eqref{eq:Ndinf} for sufficiently large $\Rc$ we have $N_\Rc \simeq C_v \Rc^d$ and hence
\begin{align*}
   \big| W_\Rc(s) - &W_\Rc(t) \big| = \frac{\Big| S_{s^\frac{1}{d} \Rc} - S_{t^\frac{1}{d} \Rc} \Big|}{\sigma_\Rc}
   \\
   &\le \frac{\Big| N_{s^\frac{1}{d} \Rc} - N_{t^\frac{1}{d} \Rc} \Big|}{\sigma_\Rc} M %
   \lesssim \frac{C_v \Rc^d \delta}{\sigma_\Rc} M \xrightarrow{\delta \to 0} 0.
\end{align*}
Therefore, $\big| W_\Rc(s) - W_\Rc(t) \big| \to 0$ as $\delta \to 0$ which,
together with~\eqref{eq:markov} ensures~\eqref{eq: tightness}.

The last property of $W_\Rc(\gamma)$ that we need to show is that of asymptotically independent increments, namely, that for all
\begin{align*}
    0 \le s_1 \le t_1 < s_2 \le t_2 < \dots < s_k \le t_k \le 1
\end{align*}
and for all linear Borel sets $H_1, \dots, H_k$ the difference
\begin{align}\label{eq: asym_ind_incr}
    \notag \P\Big( W_\Rc(t_\ell) - W_\Rc(s_\ell) &\in H_\ell, \ell = 1, \dots, k \Big) -
    \\
    &\!\!\prod_{\ell=1}^k \P \Big( W_\Rc(t_\ell) - W_\Rc(s_\ell) \in H_\ell \Big)
\end{align}
converges to $0$ as $\Rc \to \infty$. %
For any two $\sigma$-algebras $\Ac$ and $\Bc$, let
\begin{align*}
    \alpha(\Ac, \Bc) = \sup_{A \in \Ac, \ B \in \Bc} \Big | \mathbb{P} (A \cap B) - \mathbb{P}(A)\mathbb{P}(B) \Big|,
\end{align*}
and let $E_\ell$ be the event $\{ W_\Rc(t_\ell) - W_\Rc(s_\ell) \in H_\ell \}$ for any $\ell = 1, \dots, k$. Then, by definition,
\begin{align}\label{eq:Esets}
    \Big |\P\big(E_\ell \cap E_{\ell-1}\big) &- \P\big(E_\ell\big) \P\big(E_{\ell-1}\big) \Big |
    \\
    \notag &\le \alpha \bigg(\sigma \Big( \big\{\xi_i \ \big| \ s_\ell^\frac{1}{d} \Rc < ||\rbf_i|| \le t_\ell^\frac{1}{d} \Rc \big\} \Big),
    \\
    \notag &\qquad\quad \sigma \Big( \big\{\xi_j \ \big| \ s_{\ell-1}^\frac{1}{d} \Rc < ||\rbf_j|| \le t_{\ell-1}^\frac{1}{d} \Rc \big\} \Big)\bigg).
\end{align}
Furthermore, from~\cite[Prop 3.11]{bradley2007introduction}, for any two $\sigma$-algebras $\Ac$ and $\Bc$,
\begin{align}\label{eq:alpha-rho}
    \alpha(\Ac, \Bc) \le \frac{1}{4} \rho(\Ac, \Bc).
\end{align}
Therefore, combining~\eqref{eq:Esets}, \eqref{eq:alpha-rho}, and~\eqref{eq: rho_d_to_zero},
\begin{align*}
    & \Big |\P\Big(E_\ell \cap E_{\ell-1}\Big) - \P\Big(E_\ell\Big)\mathbb{P}\Big(E_{\ell-1}\Big) \Big|
    \le \frac{1}{4} \varrho^*(r_\ell)
\end{align*}
where $r_\ell = \Rc \big(s_\ell^\frac{1}{d} - t_{\ell-1}^\frac{1}{d}\big)$.
Therefore, as $\Rc \to \infty$ so does $r_\ell$, and hence
\begin{align*}
    \Big | \P\Big(E_\ell \cap E_{\ell-1}\Big) - \P\Big(E_\ell\Big)\P\Big(E_{\ell-1}\Big) \Big| \to 0 \quad \text{as} \quad \Rc \to \infty.
\end{align*}
Then, by induction,
\begin{align*}
    \bigg | \P\Big(\bigcap_{\ell=1}^k E_\ell \Big) - \prod_{\ell=1}^k \P\Big(E_\ell\Big) \bigg | \to 0 \quad \text{as} \quad \Rc \to \infty,
\end{align*}
which is the same as~\eqref{eq: asym_ind_incr}.

Put together, it follows from~\eqref{eq:mean-3d-0}, \eqref{eq:var-gamma}, \eqref{eq:uniint}, \eqref{eq: tightness}, and~\eqref{eq: asym_ind_incr} that $W_\Rc(\gamma)$ satisfies all the assumptions of~\cite[Thm 19.2]{billingsley1968convergence}, and hence
\begin{align*}
    W_\Rc(\gamma) \stackrel{d}{\to} W(\gamma) \quad \text{as} \quad \Rc \to \infty.
\end{align*}
Then, setting $\lambda = 1$ we get
\begin{align*}
    W_\Rc(1) = \frac{S_{\Rc}}{\sigma_\Rc} \stackrel{d}{\to} W(1) = \Nc(0, 1) \quad \text{as} \quad \Rc \to \infty,
\end{align*}
which completes the proof.
\end{proof}

A minor difference between Theorem~\ref{thm:clt_rho_star} and Proposition~\ref{prop:clt} is the difference between assumptions~\eqref{eq:unif_int_3d} and~\eqref{prop:cond_3}. The latter are both technical assumptions needed to prove uniform integrability of families of partial sum processes such as that in~\eqref{eq:uniint}, which in turn prevents pathological cases where some probability mass gradually drifts to infinity as $N$ or $\Rc \to \infty$. Both assumptions are sufficient but not necessary for uniform integrability, neither is uniformly stronger or weaker than the other, and either can be replaced by the direct (but perhaps more opaque) assumption of uniform integrability.

Similar to Section~\ref{sec:s_avg_1d}, we next generalize Theorem~\ref{thm:clt_rho_star} to multivariate sequences. The proof technique is the same as that of Theorem~\ref{thm:multivar_clt_rho_mixing} and hence omitted.

\begin{theorem}\longthmtitle{Multivariate CLT for $\rho^*$-mixing sequence}\label{thm:multivar_clt_rho_star_mixing}
Consider a $\rho^*$-mixing sequence of random vectors $\xib_1, \xib_2, \dots \in \real^q$ with corresponding Euclidean positions $\rbf_1, \rbf_2, \dots \in \real^d$ satisfying~\eqref{eq:Ndinf}. For any $\Rc > 0$ let $S_\Rc = \sum_{i: \|\rbf_i\| \le \Rc} \xib_i$, $N_\Rc$ be as in~\eqref{eq:NR}, $h(N_\Rc)$ be a slowly-varying residual factor for $\{\xib_i\}_{i = 1}^\infty$, and $\Sigmab_{\bar \xi}^* = \lim_{R \to \infty} \cov(S_\Rc) / N_\Rc h(N_\Rc)$. Assume, further, that
\begin{align*}
    &\mathbb{E}[\xib_i] = 0 \ \text{and} \ \|\xib_i\| \le M
    \ \text{for all} \ i \ \text{and some} \ M < \infty, \qquad\qquad
    \\
    &\sup_{\Rc > 0} \frac{\E\big[|\thetab^T S_{\gamma \Rc}|^p\big]}{\big[\thetab^T \cov(S_\Rc) \thetab\big]^\frac{p}{2}}  < \infty \ \text{for all} \ \thetab \in \real^q, \gamma \in [0, 1],
\end{align*}
and some $p > 2$.
Then,
\begin{align*}
    \bar \xib = \frac{1}{\sqrt{N_\Rc h(N_\Rc)}} \sum_{i: \|\rbf_i\| \le \Rc} \xib_i
    \stackrel{d}{\longrightarrow} \Nc(\zeros, \Sigmab_{\bar \xi}^*) %
\end{align*}
as $\Rc \to \infty$. \oprocend
\end{theorem}

\subsection{Asymptotic Linearity Under Spatial Averaging}

Theorem~\ref{thm:multivar_clt_rho_star_mixing} provides the basis that we need for solving Problem~\ref{prob2}. However, note that we have so far assumed, for simplicity, that the Euclidean locations $\{\rbf_i\}_{i = 1}^\infty$ are deterministic and fixed--e.g., at the vertices of a lattice--whereas in reality spatially-embedded subsystems are often distributed randomly over space. Therefore, before getting to the main result of this section, we prove the following lemma showing that our only requirement on the subsystems' locations, i.e., assumption~\eqref{eq:Ndinf}, is still satisfied with probability 1 in a proper stochastic setting.

\begin{lemma}\longthmtitle{Growth rate of number of subsystems distributed according to Poisson point process}\label{lem:poisson}
Consider a sequence of dynamical subsystems the locations $\{\rbf_i\}_{i = 1}^\infty$ of which are distributed according to a homogeneous Poisson point process with rate $\lambda$ on $\mathbb{R}^d$. For any $\Rc > 0$, let (the random variable) $N_\Rc$ be as in~\eqref{eq:NR}. Then
\begin{align}\label{eq:poisson}
\P \Big(\lim_{\Rc \to \infty} \frac{N_\Rc}{\Rc^d} = \lambda C_v \Big) = 1,
\end{align}
where $C_v = \frac{\pi^\frac{d}{2}}{\Gamma(1+\frac{d}{2})}$ is the volume of the unit ball in $\real^d$.
\end{lemma}
\begin{proof}
For simplicity of notation, let $V_\Rc = C_v \Rc^d$ be the volume of a ball with radius $\Rc$, and for any $r > 0$ let
\begin{align*}
\Rc(r) = \Big(\frac{r}{C_v}\Big)^\frac{1}{d}.
\end{align*}
Clearly $V_{\Rc(r)} = r$.
For all $0 \le r_0 < r_1 < \cdots < r_k$, define the annuli
\begin{align*}
A_\ell = \big\{\rbf \in \real^d \ \big| \ R(r_{\ell-1}) < \|\rbf\| \le R(r_\ell)\big\}, %
\end{align*}
and let $|A_\ell|$ be the volume of $A_\ell$. Because $A_\ell$'s are disjoint, by definition~\cite[Def. 3.1]{last2018lectures},
\begin{align*}
N_{\Rc(r_\ell)} - N_{\Rc(r_{\ell-1})}, \qquad \ell = 1, \dots, k
\end{align*}
are mutually independent Poisson-distributed random variables with rates
\begin{align*}
\lambda |A_\ell| = \lambda \big(V_{R(r_\ell)} - V_{R(r_{\ell-1})}\big) = \lambda \big(r_\ell - r_{\ell -1}\big).
\end{align*}
Therefore, $N_{\Rc(r)}, r > 0$
is a standard Poisson point process with rate $\lambda$ on the positive real line.
Let $\delta_1, \delta_2, \dots$ be the i.i.d. inter-arrival times of this standard Poisson process, and $t_\ell = \delta_1 + \cdots + \delta_\ell$. Each $\delta_\ell$ is exponentially distributed with rate $\lambda$, and hence by the strong law of large numbers,
\begin{align}\label{eq:slln}
\frac{t_\ell}{\ell} = \frac{\delta_1 + \cdots + \delta_\ell}{\ell} \stackrel{a.s.}{\to} \frac{1}{\lambda} \quad \text{as} \quad \ell \to \infty.
\end{align}
For any $r \in [t_\ell, t_{\ell+1})$, by definition $N_{\Rc(r)} = \ell$ and thus
\begin{align}\label{eq:upper-lower}
\frac{\ell}{t_{\ell+1}} \le \frac{N_{\Rc(r)}}{r} \le \frac{\ell}{t_\ell}.
\end{align}
By~\eqref{eq:slln}, as $r \to \infty$ (and thus $\ell \to \infty$) both sides of~\eqref{eq:upper-lower} converge to $\lambda$ a.s., and thus
\begin{align}\label{eq:aslimit}
\frac{N_{\Rc(r)}}{r} \stackrel{a.s.}{\to} \lambda \quad \text{as} \quad r \to \infty.
\end{align}
The statement of the theorem then follows by replacing $r = C_v \Rc(r)^d$ in~\eqref{eq:aslimit} and changing the limiting variable from $r$ to $\Rc$.
\end{proof}

Lemma~\ref{lem:poisson} allows us to restate the preceding CLTs for $\rho^*$-mixing sequences with fixed Euclidean locations for the case where the locations are Poisson distributed. %
These restatements are straightforward and hence omitted. Instead, we proceed to use the combination of Lemma~\ref{lem:poisson} and the preceding CLTs to solve Problem~\ref{prob2}. The following is the parallel to Assumption~\ref{assum} for spatially-embedded systems.

\begin{assumption}\longthmtitle{Updated assumptions for spatially-embedded systems}\label{assum-star}
Consider the population of dynamical systems described by~\eqref{eq:ode} and~\eqref{eq:x_Ni_3d} and let $\ybf_i(t)$, $\zbf_i(t)$, and $N_\Rc$ be as in~\eqref{eq:rho-mixing-1d}, \eqref{eq: conc_states}, and~\eqref{eq:NR}, respectively. We make the following assumptions.
\begin{enumerate}[label=\rm(A\arabic*$^*$), leftmargin=25pt]
    \item\label{assum:w-star} $\E[\wbf_i(t)] = \zeros$ and $\|\wbf_i(t)\| \le C$ a.s. for all $i, t$ and some $C < \infty$.

    \item\label{assum:st-var-growth-star}
    For all $\thetab \neq 0$, $t \ge 0$, and $\gamma \in [0, 1]$ and some $p > 2$,
    \begin{align*}
    \sup_{\Rc > 0} \frac{\E\big[\big|\thetab^T \sum_{i: \|\rbf_i\| \le \gamma \Rc} \zbf_i(t)\big|^p\big]}{\big[\thetab^T \cov\big(\sum_{i: \|\rbf_i\| \le \Rc} \zbf_i(t)\big) \thetab\big]^\frac{p}{2}}  < \infty.
    \end{align*}

    \item\label{assum:w-x-indep-star} Same as Assumption~\ref{assum:w-x-indep}.

    \item\label{assum:f-bdd-star} Same as Assumption~\ref{assum:f-bdd}.

    \item\label{assum:rhomixing-star} For all $t$, the sequence of random vectors $\{\ybf_i(t)\}_{i = 1}^\infty$ is $\rho^*$-mixing with a slowly-varying residual factor $h(N_\Rc)$. \oprocend
\end{enumerate}
\end{assumption}

We are now ready to state the main result of this section, as follows. The proof is similar to the proof of Theorem~\ref{thm:s_avg} and hence omitted.

\begin{theorem}\longthmtitle{Linearizing effect of spatial averaging on spatially-embedded populations of dynamical systems}\label{thm:s_avg_3d}
Consider the population of dynamical systems described by~\eqref{eq:ode} and~\eqref{eq:x_Ni_3d}, and assume that Assumptions~\ref{assum:w-star}-\ref{assum:rhomixing-star} as well as~\eqref{eq:Ndinf} hold.
Define the population average variables $\bar \xbf(t)$ and $\bar \wbf(t)$ as in~\eqref{eq:s_avg_state_3d}, let $\bar \ybf(t) = \begin{bmatrix} \bar \xbf(t)^T & \bar \wbf(t)^T \end{bmatrix}^T$, and assume that $\lim_{\Rc \to \infty} \cov\big(\bar \xbf(t+1), \bar \ybf(t)\big)$ exists.
Then, for all $\xib \in \real^n, \omegab \in \real^m$,
\begin{align}\label{eq:cond-exp-s_avg-star}
    \!\!\!\E\bigg [ \mathbf{\bar x}(t+1) \bigg | \!\begin{bmatrix}
    \mathbf{\bar x}(t) \\ \mathbf{\bar w}(t)
    \end{bmatrix} = \begin{bmatrix}
    \xib \\ \omegab
    \end{bmatrix} \!\bigg ]
    \!\xrightarrow{\Rc \to \infty}\! \Abf_\infty(t)\xib  + \Bbf_\infty(t) \omegab
\end{align}
where $\Abf_\infty(t)$ and $\Bbf_\infty(t)$ are defined in~\eqref{eq:s_AB}. Furthermore, if the subsystem locations $\{\rbf_i\}_{i = 1}^\infty$ are distributed according to a homogeneous Poisson point process on $\mathbb{R}^d$, \eqref{eq:cond-exp-s_avg-star} holds with probability 1. \oprocend
\end{theorem}

As we noted in the Introduction,the present work was motivated by our empirical observations of macroscopic linearity in brain networks. Theorem~\ref{thm:s_avg_3d} now rigorously explains these observations, and highlights the key roles of spatial averaging that is inherent in macroscopic recordings such as electroencephalography (EEG)~\cite{buzsaki2012origin} and functional magnetic resonance imaging (fMRI)~\cite{buxton2013physics}, as well as the weak and decaying correlations that have long been observed between biological neurons~\cite{smith2008spatial,smith2013spatial,rosenbaum2017spatial}.
Whether and to what extent large-scale complex networks in other domains meet the assumptions of this theorem, and are hence bound to follow its consequent macroscopic linearity, remains an interesting topic of inquiry for future research.

\section{Conclusions}\label{sec:conclusions}

In this article we developed a theoretical framework to understand the linearity of spatially averaged dynamics in heterogeneous populations of nonlinear networked dynamical systems. This research was motivated by several observations of linear behavior at the macroscopic level in biological and artificial neural systems, as well as in-silico observations that nonlinear systems nearly universally exhibit (more) linear dynamics when subjected to spatial averaging. To our knowledge, this work is the first to formalize these observations into a unified theoretical framework and support them with rigorous mathematical analysis. %
By building on and extending the celebrated central limit theorem and the theory of mixing sequences, we proved that averaging has a strong and robust linearizing effect that holds for almost any form of microscopic nonlinearity, noise distribution, and network connectivity patterns, \textit{as long as} pairwise correlations decay with some notion of distance between microscopic subsystems.
We proved this result in two general settings, one where pairwise correlations decay with the distance $|i-j|$ between subsystems' linear indices, and the other when pairwise correlations decay with the Euclidean distance $\|\rbf_i - \rbf_j\|$ between the physical locations of spatially-embedded subsystems. These results were further extended to time-invariant limit dynamics, finite-sample averaging with rates of convergence, and networks of spatially-embedded subsystems with random locations.
Overall, our results provide significant insights into the macroscopic behavior of large-scale systems, and lay a robust theoretical foundation for future research across various domains. %

\begin{ack}
This work was support in part by the National Science Foundation Award No. 2239654 to EN.
\end{ack}

\end{document}